\documentclass[12pt]{amsart}

\def\@typesizes{%
       \or{5}{6.5}\or{6}{7.5}\or{7}{8.5}\or{8}{11}\or{9}{12}%
       \or{10}{13}
       \or{\@xipt}{14}\or{\@xiipt}{15}\or{\@xivpt}{18}%
       \or{\@xviipt}{20}\or{\@xxpt}{24}}

\usepackage{amsthm}
\usepackage{amsmath}
\usepackage{amssymb}
\usepackage{graphicx}
\usepackage{enumerate}
\usepackage{tikz}
\allowdisplaybreaks
\numberwithin{equation}{section}
\numberwithin{figure}{section}

\theoremstyle{plain}
\newtheorem{theorem}{ Theorem}[section]
\newtheorem{proposition}[theorem]{ Proposition}
\newtheorem{lemma}[theorem]{ Lemma}
\newtheorem{corollary}[theorem]{ Corollary}
\newtheorem{example}[theorem]{ Example}
\newtheorem{remark}[theorem]{ Remark}
\newtheorem{definition}[theorem]{ Definition}
\newtheorem{conjecture}{ Conjecture}

\def\BET{\begin{theorem}}
\def\ENT{\end{theorem}}
\def\BEP{\begin{proposition}}
\def\ENP{\end{proposition}}
\def\BEL{\begin{lemma}}
\def\ENL{\end{lemma}}
\def\BEC{\begin{corollary}}
\def\ENC{\end{corollary}}
\def\BEE{\begin{example} \rm}
\def\ENE{\end{example}}
\def\BER{\begin{remark} \rm}
\def\ENR{\end{remark}}
\def\BED{\begin{definition} \rm}
\def\END{\end{definition}}
\def\BECJ{\begin{conjecture}}
\def\ENCJ{\end{conjecture}}

\def\bea{\begin{eqnarray}}
\def\eea{\end{eqnarray}}

\def\beas{\begin{eqnarray*}}
\def\eeas{\end{eqnarray*}}

\def\beq{\begin{equation}}
\def\eeq{\end{equation}}

\def\beal{\begin{align*}}

\def\eeal{ \end{align*} }

\def\roweq{\nonumber \\ &=& }
\def\rowleq{\nonumber \\  & \leq & }
\def\rowgeq{\nonumber \\ & \geq & }

\def\rowpl{\nonumber \\ &  \ \ + & }

\def\bfS{{\bf S}}

\def\bbC{{\mathbb C}}
\def\bbD{{\mathbb D}}

\def\bbN{{\mathbb N}}

\def\bbR{{\mathbb R}}

\def\bbZ{{\mathbb Z}}

\def\cH{{\mathcal H}}
\def\cI{{\mathcal I}}

\def\cL{{\mathcal L}}
\def\cM{{\mathcal M}}

\def\cP{{\mathcal P}}

\def\cR{{\mathcal R}}

\def\cT{{\mathcal T}}

\def\cX{{\mathcal X}}

\def\sF{{\sf F}}
\def\wP{{\widetilde P_\mu}}
\def\ef{\eqref}

\def\sfa{\mathsf{a}}
\def\sfb{\mathsf{b}}
\def\sfh{\mathsf{h}}
\def\sft{\mathsf{t}}

\begin{document}

\title[Toeplitz operators in periodic domains]{On Bergman-Toeplitz operators in periodic planar domains}

\author{Jari Taskinen}

\address{Department of Mathematics and Statistics,
 University of Helsinki, P.O. Box 68, 00014 Helsinki, Finland}
 
\email{jari.taskinen@helsinki.fi}

\thanks{ This research was partially supported by the V\"ais\"al\"a Foundation of
the Finnish Academy of Sciences and Letters and by the Academy of Finland
project no. 342957. The author also wishes to thank Antti Per\"al\"a (Ume{\aa} 
University) for some discussions concerning the example in Section \ref{sec9R}.}

\subjclass{Primary 47A10, 47B35; secondary 30H20, 30C40, 47A75, 47B91}

\begin{abstract}
We study spectra of Toeplitz operators $T_a $ with periodic
symbols in Bergman spaces $A^2(\Pi)$ on  unbounded periodic  planar 
domains $\Pi$, which are defined as the union of infinitely many copies of
the translated, bounded periodic cell $\varpi$. We introduce 
Floquet-transform techniques and prove a version of the band-gap-spectrum 
formula, which
is well-known in the framework of periodic elliptic spectral problems and 
which describes the essential spectrum of $T_a$ in terms of the spectra of a 
family of Toepliz-type operators $T_{a,\eta}$ in the cell $\varpi$, where
$\eta$ is the so-called Floquet variable. 

As an application, we consider periodic domains $\Pi_h$ containing thin
geometric structures and show how to construct a Toeplitz operator $T_\sfa: 
A^2(\Pi_h) \to A^2(\Pi_h)$ such that the essential  spectrum of $T_\sfa$ 
contains disjoint components  which approximatively coincide with any given 
finite set of real numbers. Moreover, our method provides a systematic and 
illustrative way how to construct such examples by using Toeplitz operators
on the unit disc $\bbD$ e.g. with radial symbols.  

Using a Riemann mapping one can then find a 
Toeplitz operator $T_a : A^2(\bbD) \to A^2(\bbD)$ with a bounded symbol and  with 
the same spectral properties as $T_\sfa$.

\end{abstract}

\maketitle

\section{Introduction.}
\label{sec1}

The spectral theory of  Toeplitz operators $T_a$  in Bergman spaces $A^2(\Omega)
$ even on  the unit disc $\Omega = \bbD$ is not well understood yet, especially 
concerning the essential spectrum $\sigma_{\rm ess} (T_a)$ of $T_a$. The purpose of
this paper is to apply Floquet-transform methods and prove formula \ef{0.0}
which connects the essential spectrum of a Toeplitz operator on a periodic domain
$\Pi$ and the spectra of a family of Toeplitz-type operators on the periodic cell
$\varpi$ of $\Pi$. The corresponding formula \ef{8.6} is fundamental in the 
completely different setting of the Schr\"odinger 
equation with periodic potentials, since it establishes the band-gap structure of the
spectrum; see \cite{Ku}, Sections 5.11. and 6.5. for the connections of 
the formula to Fermi surfaces and Wannier functions of molecular physics. 
The formula is also extensively used in the study of elliptic spectral problems in 
periodic domains, see for example \cite{BaTa}, \cite{CCNT}, \cite{N1}, \cite{NaPl}, 
\cite{NaTa} and many others. 

In the case of Bergman-Toeplitz-operators,  formula \ef{0.0} will  yield a 
new approach to spectra and in particular a systematic method for new types of 
examples, see Theorem \ref{cor9.2} and its proof. To recall some of the known
result in this area, in the relatively simple case that the symbol $a$ is 
harmonic on $\bbD$ and has
a continuous extension to the closed disc there holds $\sigma_{\rm ess} (T_a )
= a( \partial \bbD)$.  Moreover,  $\sigma_{\rm ess}(T_a)$  is connected, if 
$a$ is harmonic and either real-valued or piecewise continuous on $\partial 
\bbD$, see \cite{MS}. Also, in \cite{StZ}, Corollary 20, it was shown 
that $\sigma_{\rm ess}(T_a)$ is connected, if the symbol belongs to the algebra
of those $a \in L^\infty(\bbD)$ such that the corresponding Hankel
is compact, thus, for example, if the symbol belongs to the class $VMO_\partial(\bbD)$ 
(see Section 8.4. of \cite{Zh}). 
A harmonic symbol such that $\sigma_{\rm ess}(T_a)$ is disconnected was constructed in
\cite{SuZ}. In this example, 0 is an isolated point of the essential spectrum. 
In the case of the Hardy space, it is known that the essential spectrum of a bounded Toeplitz operator is always connected, see \cite{Doug}.

On the other hand, there are certain special cases where one can get a 
information on the spectrum of $T_a$ in a  straightforward way and 
at least calculate some or all eigenvalues. For example, if the symbol 
$a$ is radially symmetric, i.e. $a(z) = a(|z|)$ for all $z \in \bbD$, then $T_a$ is 
the Taylor coefficient multiplier $T_a : \sum_{n=0}^\infty f_n z^n \mapsto
\sum_{n=0}^\infty \lambda_n f_n z^n $, where the eigenvalues are 
\bea
\lambda_n = \frac{n+1}{\pi}\int\limits_0^1 a(r) r^{2n+1} dr.  \label{Taylor}
\eea 
To describe the spectrum there remains to characterize the numbers 
$\lambda_n$, which, in spite of the concrete moment formula, may be a 
complicated task,  depending on what is wanted. See \cite{V1}, \cite{V2}, 
\cite{V3}. 
Our considerations of Toeplitz operators on periodic planar domains 
$\Pi \subset \bbC$ will allow us to transform some properties of the spectra of 
these type of operators and others on $\bbD$ into results on the
essential spectra of periodic symbols. By applying the Riemann conformal mapping
from $\Pi$ onto $\bbD$ (Lemma \ref{lem1.5}) one obtains new types of interesting 
examples of Toeplitz operators in   $A^2( \bbD)$. These results will be consequences of a 
systematic treatment of Toeplitz-operators  in the periodic setting via Floquet-transform 
methods, based on \cite{T1} and  on the corresponding machinery in elliptic spectral 
problems, \cite{Kbook}, \cite{Ku}, \cite{N1}, \cite{NaPl}.

The {\it unbounded} periodic domain $\Pi$ is obtained as the union of infinitely  many translated copies of the {\it bounded} periodic 
cell $\varpi \subset \bbC$ (more precisely, as the interior of the union of the closures 
of the translates). We will consider periodic symbols $a$ on $\Pi$ and, as
mentioned,  the main aim is to 
establish the following connection of the spectrum $\sigma (T_a)$ and 
the essential spectrum $\sigma_{\rm ess} (T_a)$ of 
$T_a$ on  $A^p(\Pi)$ with the spectra $\sigma (T_{a,\eta})$ of a family of related 
Toeplitz operators $T_{a,\eta}$ on Bergman-type spaces $ A_\eta^2(\varpi)$ defined on 
the periodic cell, see Theorem \ref{th8.3}:
\bea
\sigma( T_a )=\sigma_{\rm ess} (T_a) = \bigcup_{\eta \in [-\pi,\pi]} \sigma (T_{a,\eta}) .  
\label{0.0}
\eea
Here,  $\eta \in [-\pi,\pi]$ is the so-called Floquet-parameter. The 
corresponding formula for elliptic spectral problems in periodic domains can
be found in Theorem 2.1. of \cite{N1} and Theorem 3.4.6. of \cite{NaPl} and in the
case of the Schr\"odinger operator with periodic potentials and other elliptic 
operators with periodic coefficients in  \cite{Kbook},  and \cite{Ku},  Theorem 5.5. 
These references provide some guidance to the topic, but proving the 
result and working in spaces of analytic functions is certainly a world different from
pde-arguments in Sobolev-spaces of functions of real variables. 
A second main source is the paper
\cite{T1}, where the author established the basic theory of the Floquet-transform and
Bergman projection  on periodic domains.  

The second main result, Theorem \ref{cor9.2}, is an application of the previous 
considerations and it states that given a finite sequence of real
numbers $(\lambda_n)_{n=1}^N$, there exists a Toeplitz-operator $T_a$ on the Bergman space 
of the unit disc $\bbD$ such that its essential spectrum $\sigma_{\rm ess} (T_a)$ consists 
of arbitrarily close approximations of the numbers $\lambda_n$ and other components, 
which are at some distance of these numbers. Thus, $\sigma_{\rm ess} (T_a)$ 
can for example be made to have arbitrarily many disjoint components.\footnote{We expect
these components to be continua, but we do not have a proof for this. The corresponding
question in elliptic pde-theory is quite deep and partially open, see Remark $2^\circ$ below
Theorem \ref{th9.1}.} 
The construction
of such operators is first made on a periodic domain $\Pi$, where the periodic 
cells are connected by ligaments of width $2h$ with a small $h > 0$, and then 
using formula \ef{0.0}. 
Finally, applying a Riemann mapping  allows us to transform
the operator to the unit disc. 

Due to technical reasons, which for example include the need to stick to self-adjoint 
operators, our concrete examples will concern real-valued symbols,
the restrictions of which to the periodic cells are compactly supported. Thus, these symbols
are far from being harmonic. On the other hand, the general Theorem \ref{th8.3} 
does not include restrictions as regards to harmonicity so that it will remain 
open if examples like in Theorem \ref{cor9.2} could be constructed for harmonic 
symbols. 

Theorem \ref{cor9.2} also has its predecessors in the theory of elliptic spectral problems
related to unbounded self-adjoint operators in Hilbert-spaces, nam\-ely,
there are many works where the existence of gaps in the essential spectrum are studied
in periodic domains with thin ligaments, for example see \cite{BaTa}, \cite{CCNT}, \cite{NaRuTa}, \cite{NaTa}.
We borrow from these papers the idea that 
the essential spectrum of the problem in the full, periodic domain $\Pi$ is 
an approximation of the spectrum (eigenvalues) of a problem in the bounded periodic cell
$\varpi$: formula \ef{0.0} is one of the key ingredients here, but there are others, and 
again, there of course are major methodological differences between elliptic 
pde-considerations and arguments in Bergman spaces.

Let us proceed to some basic definitions, notation and preliminary results for this paper.  
Given a domain $\Omega$ in the complex plane $\bbC$, we
denote by $L^2(\Omega)$ the usual Lebesgue-Hilbert space with respect to 
the (real) area measure $dA$ and by  $A^2(\Omega)$ the corresponding Bergman space, 
which is the subspace consisting of analytic functions. The norm of $f \in L^2(\Omega)$ 
is denoted by  $\Vert f \Vert_\Omega$ and the inner product of $f,g\in L^2(\Omega)$ by 
$(f | g )_\Omega = \int_\Omega f \bar g$. 
It is a consequence of the
Cauchy integral formula that the norm topology of $A^2(\Omega)$ is stronger 
than the topology of the uniform convergence on compact subsets, and this 
implies that the Bergman space is always a closed subspace of $L^2(\Omega)$,
hence it is a Hilbert space, in particular complete. We denote by $P_\Omega$
the orthogonal projection from $L^2(\Omega)$ onto $A^2(\Omega)$. It can 
always written with the help of the Bergman kernel $K_\Omega :
\Omega \times \Omega \to \bbC$,
\beas
P_\Omega f(z) = \int\limits_\Omega K_\Omega (z,w) f(w) dA(w)  
\eeas
and the kernel has the properties that $K_\Omega (z, \cdot) \in L^2(\Omega)$ 
for all $z$ and $K(z,w) = \overline{ K(w,z)}$ for all $z,w \in \Omega$. 
See e.g. \cite{K1} for a proof of these assertions. 

Given a function $a \in L^\infty (\Pi)$, the Toeplitz operator $T_a : A^2(\Pi) \to A^2(\Pi)$  
with symbol $a$ is defined by
\bea
T_a f(z) = P_\Omega M_a f(z) = \int\limits_\Omega K_\Omega (z,w) a(w) f(w) dA(w)  
\label{1.1a}
\eea
where $M_a$ is the pointwise multiplication operator $ f \mapsto af$.  It is 
plain that $T_a$ is a well-defined bounded linear operator, since $M_af \in L^2(\Pi)$ 
by  our assumptions.

The following notation will be used throughout the paper. We write $C$, $C'$, $\ldots$,
(respectively, $C_\eta$, $C'_\eta$, $\ldots$ etc.) for positive constants independent  
of  functions or variables in the given inequalities (resp. depending only on a 
parameter  $\eta$ etc.), the values of which may vary from place 
to place. For $x \in \bbR$, $[x]$ denotes the largest integer not larger than
$x$. If $x$ is a point in $ \bbC$ or $\bbR^n$ and $r >0$, then $B(x,r)$ 
denotes the Euclidean ball with center $x$ and radius $r>0$. Moreover,
cl$(A)$ denotes the closure of a set $A$. Given a function $f$, its support is
denoted by supp\,$f$. We write $\bbZ= \{ 0, \pm1, \pm2, \ldots \}$. 
Given $R$ with $0 < R < 1$ we denote $ \{ z \in \bbC \, : \, |z| < R \} =: \bbD_R$.

We use standard notation and definitions for $L^p$-spaces with $1\leq p \leq \infty$,
in particular  $L^\infty(\Omega)$ stands for the Banach-space of essentially bounded complex 
valued  functions on  $\Omega$.
In general, the norm of a Banach-space $X$ is denoted by $\Vert \cdot \Vert_X$. 
Given an interval $I \subset \bbR$ and a Banach space $X$ we denote  by $L^2(I;X)$ the 
space of vector valued, Bochner-$L^2$-integrable functions on $I$ with values  in  $X$, 
endowed  with the norm 
\beas
\Vert f \Vert_{L^2(I;X)}=\Big(\int\limits_I 
\|f( t )\|_X d t \Big)^{1/2}\, .
\eeas 
If $X = L^2(\Omega)$ for some domain $\Omega$, then $L^2(I;X)$ is a Hilbert space 
endowed with the  inner product $\int_I (f (t)| g (t))_\Omega dt  $.
See \cite{Hyt} for the theory of Bochner-spaces.

Given a Banach space $X$, $\cL(X)$ stands for the Banach space of bounded
linear operators $X \to X$. The operator norm of $T \in \cL(X)$ is  denoted just by
$\Vert T \Vert$ or by $\Vert T \Vert_{X \to X}$, if it is necessary to specify the 
domain or target spaces. The identity operator on $X$ is written as $I_X$ or by
$I_\Omega$, if $X=L^2(\Omega)$.

For an operator  $T \in \cL(H)$, where $H$ is a Hilbert-space,  $\sigma (T)$, $\sigma_{\rm ess} (T)$ and 
$\varrho(T)$ stand for the spectrum, essential spectrum and resolvent set of the 
operator $T$. The resolvent (operator) of  $T \in \cL(H)$ is denoted by 
$R_\lambda(T) = (T -\lambda I_H)^{-1}$, where $\lambda \in \varrho(T)$.

We finish this section by recalling known facts which will be needed in 
the proof of Theorem \ref{cor9.2}. The first one  is the Weyl criterion 
for the points of the essential spectrum, see e.g. Theorem VII.12 in \cite{RS}. 

\BEL \label{lemW}
If $H$ is Hilbert space, $T\in \cL(H)$ and $\lambda \in \bbC$, then $\lambda \in 
\sigma_{\rm ess}(T)$, if and only if there exists a Weyl singular sequence, 
which is a sequence $(h_n)_{n=1}^\infty$ of elements of $H$ such that $\Vert h_n \Vert_H =1$
for all $n$, 
\bea
\lim\limits_{n\to \infty} \Vert T h_n - \lambda h_n\Vert_H = 0
\label{8.39}
\eea
and such that the sequence has no convergent subsequences. 
\ENL

The next is the lemma on almost eigenvalues and -vectors, which is a consequence
of the spectral theorem. See \cite{VL}, or  also for example Lemma 4.2. in \cite{BaTa}, 
or Lemma 5.3. in \cite{CCNT}. 

\BEL
\label{lem9.2}
Let $K: H \to H$ be a compact self-adjoint operator in a Hilbert-space $H$ and 
let $\mu \in \bbR$. If there are  $f \in H$ with $\Vert  f\Vert_H =1$ and  $\delta > 0 $ 
such that 
\beas
\Vert K f - \mu f \Vert \leq \delta  , 
\eeas
then $K$ has an eigenvalue $\lambda \in [\mu - \delta, \mu + \delta]$. 
\ENL

The kernel formula appearing in \ef{1.12} of the following lemma can be found e.g. in 
\cite{Be}, formula (1), or \cite{K1}, Proposition 2.7.  All 
assertions of Lemma \ref{lem1.5} are well-known and can be easily checked by standard calculations.

\BEL
\label{lem1.5}
Let $\Omega \subset \bbC$ be a simply connected domain and let $\varphi $
be the Riemann conformal mapping from $ \Omega$ onto $\bbD$ with inverse
$\psi = \varphi^{-1} : \bbD \to \Omega$. Then, the operator
$L: f \mapsto \psi' f \circ\psi$ is an isometric Hilbert-space isomorphism
from the Bergman space $A^2 (\Omega) $ onto $A^2(\bbD)$ and
$L^{-1} =: J $ is the mapping $f \mapsto \varphi' f \circ \varphi$. 

Moreover, if $\sfa \in L^\infty(\Omega)$ and $T_\sfa: A^2(\Omega) \to A^2(\Omega)$ is the
(bounded) Toeplitz operator with symbol $\sfa$, i.e.
\bea
T_\sfa f(z) &=& \int\limits_\Omega K_\Omega(z,w)\sfa(w) f(w) dA(w) 
\roweq 
\int\limits_\Omega \frac{\varphi'(z) \overline{\varphi'(w)} \sfa(w)f(w) }{(1- \varphi(z) 
\overline{ \varphi(w)})^2} dA(w) ,   \label{1.12}
\eea
then $L T_\sfa J$ is the Toeplitz operator $T_a : A^2(\bbD) \to A^2(\bbD)$ with the bounded 
symbol $a = \sfa \circ \psi$. 
\ENL

\BER \label{lem9.3}
Also, the following simple consequence of the Pythagorean theorem will be used:
if  $X \subsetneq H$ is  a closed subspace of a Hilbert space $H$,
$P_X$ is the orthogonal projection from $H$ onto $X$, and $g_1, g_2 
\in X$ and  $0 \not= f \in H \setminus X$ and  $\delta > 0$
are such that $P_X f = g_1$ and $\Vert g_2 - f\Vert 
\leq \Vert g_1 -f \Vert  +\delta$, then there holds
\bea
\Vert g_1 - g_2 \Vert \leq 2 \delta \Vert g_1 -f \Vert + \delta^2. 
\label{9.3a}
\eea
Indeed, since $P_X f = g_1$ and $g_2 \in X$, the vectors $g_1 - g_2$
and $g_1 - f$ are orthogonal to each other, by the definition of the orthogonal projection, and we thus  have 
\beas
\Vert g_1 - g_2 \Vert^2 + \Vert g_1 - f \Vert^2
= \Vert  g_2 - f \Vert^2 \leq \Vert g_1 -f \Vert^2  + 2 \delta
\Vert g_1 -f \Vert + \delta^2,
\eeas
which implies \ef{9.3a}.

\ENR

\section{Preliminaries on periodic domains, Floquet-transform and Bergman projection.}
\label{sec2}

In this section we recall the basic theory of Bergman spaces and projections in 
periodic planar domains $\Pi$ as presented in the paper \cite{T1}. In particular, 
the general geometric assumptions on $\Pi$ will be the same as \cite{T1}. 
Later, in the example of Sections \ref{sec9R}--\ref{secY} we will consider a
family of  such domains depending on a small geometric parameter $h>0$ (the width 
of ligaments connecting the periodic cells). 

We denote the periodic cell by $\varpi $ and require that  $\varpi \subset ]0,1[ 
\times ]-M, M[ \subset \bbC$  
for some $M>0$ and that  the boundary $\partial \varpi$ of $\varpi$ is a Lipschitz 
curve such that its  intersection  
with the axis $\{ z \, : \, \mbox{Re}\, z = \frac12 \pm \frac12 \}$
coincides with $\{ \frac12 \pm \frac12  \} \times [a,b] =: J_{\pm } $ 
for some real numbers $b > a$. By assumption, $\varpi$ is a 
Carath\'eodory domain 
so that by \cite{Ca}, polynomials form a 
dense subspace of the Bergman space $A^2(\varpi)$ (see also the 
introduction in the paper \cite{Br}). 
We denote the translates of $\varpi$ by $\varpi(m) = \varpi + m$, where $m 
\in \bbZ  \subset \bbC$, and then define the periodic domain $\Pi$ 
as the interior of the set
\bea
\bigcup_{m \in \bbZ} {\rm cl}( \varpi(m) ).  \label{1.2}
\eea
Note that $\Pi$ is simply connected, if and only $\varpi$ is. 

The Floquet-transform is defined for $f \in A^2(\Pi)$ by
\bea 
\sF f(z,\eta)&=& 
\frac{1}{\sqrt{2\pi}}\sum_{m\in \mathbb{Z}}
e^{-i \eta m  } f (z + m ), \ \ z \in \varpi, \ \eta \in [-\pi, \pi]. \label{1.14}
\eea
The series in \eqref{1.14} converges in $ L^2(-\pi,\pi; L^2(\varpi)) :=  L^2\big( [-
\pi,\pi] ; L^2(\varpi)  \big)$, thus also pointwise 
for a.e. $\eta,z$, and also in $L^2(\varpi)$ for a.e. $\eta$. In \cite{T1} there
is a simple example which shows that the series does not need to converge
for exactly all $z$. 
If  $g \in L^2(-\pi, \pi ; 
L^2(\varpi)$, we also denote
\bea
\sF^{-1} g (z ) & = &  \frac1{\sqrt{2 \pi}} \int\limits_{-\pi}^\pi
e^{ i [{\rm Re} z] \eta} g(z -  [{\rm Re}\,z ] , \eta ) d \eta ,  \ \ z \in \Pi  .
\label{1.14b}  
\eea 

The following is one of the main results of the paper \cite{T1}.

\BET
\label{prop1.2}
Floquet transform $\sF$ maps $A^2(\Pi)$ onto $ L^2(-\pi,\pi; 
A_\eta^2(\varpi))$. Its inverse 
$\sF^{-1} :  L^2(-\pi,\pi; A_\eta^2(\varpi)) \to A^2(\Pi)$ 
is given by the formula \eqref{1.14b}. Moreover, 
$\sF$ preserves the inner product and is thus a unitary operator. 
\ENT

Here, the space $A_\eta^2(\varpi)$ is defined as follows. First,
given  $\eta \in [-\pi, \pi]$, we denote by $A_{\eta, {\rm ext}}^2(\varpi)  $ 
the subspace of $ A^2(\varpi)$ consisting of functions  $f $ which can be 
extended as analytic functions to a  neighborhood in the domain $\Pi$ of the 
set ${\rm cl} \,(\varpi) \cap \Pi$    and satisfy  the  boundary condition 
\begin{equation}
f \Big( \frac12 + iy \Big) = e^{i \eta } f  \Big( -\frac12 + iy \Big)
\ \ \mbox{for all} \ a < y < b. \ 
\label{1.15}
\end{equation}
We define the space  $A_\eta^2 (\varpi) $ as the 
closure of $A_{\eta, {\rm ext}}^2(\varpi)  $ in $A(\varpi)  $. 
Note that in general, the functions belonging to 
$A^2(\varpi)$  may not have properly defined boundary values on $\partial 
\varpi$ so that the  condition \eqref{1.15} cannot be posed directly.
Finally, $ L^2(-\pi,\pi; A_\eta^2(\varpi))$ denotes the subspace of the
Hilbert-space  $ L^2(-\pi,\pi; A^2(\varpi))$ consisting of functions $f$ such that the 
function $z \mapsto f (z, \eta)$ belongs to $A_\eta^2(\varpi)$ for a.e. $\eta \in [-\pi, 
\pi]$. One observes that $ L^2(-\pi,\pi; A_\eta^2(\varpi))$ is a closed  subspace of 
$ L^2(-\pi,\pi; L^2(\varpi))$. Indeed, this will suffice for our considerations, but
we remark  that $ L^2(- \pi,\pi; A_\eta^2(\varpi))$
is not a vector valued $L^2$-space, since the space $ A_\eta^2(\varpi)$ depends on
$\eta$,  but it only has the structure of a Banach vector bundle, 
see for example Section 1.3 of \cite{Kbook}. 

\BED  \label{defPe}
The orthogonal projection  from $L^2 (\varpi) $ onto the subspace 
$A_\eta ^2(\varpi)$ is denoted by $P_\eta$.
\END

\BER \label{rem2.3}
One can show that the mapping  $\cP$, defined in the space $L^2(-\pi,\pi; L^2(\varpi))$ by 
taking the projection $P_\eta$ pointwise for all $\eta \in [-\pi,\pi]$, is 
actually the orthogonal projection from $L^2(-\pi,\pi; L^2(\varpi))$ onto
$L^2(-\pi,\pi; A_\eta^2(\varpi))$, and the Bergman projection $P_\Pi$ from 
$L^2(\Pi)$ onto $A^2(\Pi)$ can be written as $P_\Pi = \sF^{-1} \cP \sF$.
See Section 4 of \cite{T1} for the details. 
We will use in the sequel a similar presentation for Toeplitz operators in
$A^2(\Pi)$, see Lemma \ref{lem8.0a}, below.

For $P_\Pi$ and the translations $\tau_n (z) = z+ n$, $n \in \bbZ$, there holds
the commutation relation $P_\Pi (f \circ \tau_n) = (P_\Pi f ) \circ \tau_n$ for all $f \in 
L^2(\Pi)$ and $n \in \bbZ$. 
\ENR

\section{Preliminaries on Toeplitz operators with periodic symbols}  
\label{sec8}

In this section we consider Toeplitz operators $T_a$ with periodic symbols 
and in particular their essential spectra. To mention some background 
information, we consider  the following spectral Laplace-Neumann problem
\bea
- \Delta u (\bar x) &=& \lambda u (\bar x), \ \ \ \bar x \in \Pi,  \nonumber 
\\ 
\partial_\nu u(\bar x) &=& 0, \ \ \ \ \ \bar x \in \partial \Pi,\label{8.2}
\eea 
where $\Pi$ is a domain in $\bbR^2 \cong \bbC$ with the same geometric
assumptions as above, $u$ is a twice continuously differentiable
function of the real variable $\bar x =  (x,y)\in \Pi$, $\lambda$ is the spectral parameter and
$\partial_\nu$ is the outward  normal derivative on the boundary $\partial \Pi$ of
$\Pi$. In the present geometric situation, it is known that the essential spectrum
of the problem \ef{8.2} is non-empty. The Floquet transform \eqref{1.14} (more precisely, its standard real variable
version) turns problem \ef{8.2} into an $\eta$-dependent family of spectral Laplace problems: these are  obtained by
restricting problem \ef{8.2} to the periodic cell $\varpi$ and adding on the lateral
sides $J_\pm$ of $\partial \varpi$ the quasiperiodic boundary condition 
\bea
u(1/2, y) = e^{i \eta} u(-1/2,y), \ \ y \in [a,b], 
\eea
depending on the parameter $\eta \in [-\pi,\pi]$. If the spectrum of the latter problem
is denoted for each $\eta$ by $\sigma(\eta)$, then there holds the well-known formula
\bea
\sigma_{\rm ess} = \bigcup_{\eta \in [- \pi, \pi] } \sigma(\eta)  \label{8.6}
\eea
for the essential spectrum $\sigma_{\rm ess}$ of the problem \ef{8.2};
see \cite{Kbook},  and \cite{Ku},  Theorem 5.5., \cite{N1} , Theorem 2.1, and  \cite{NaPl}, 
Theorem 3.4.6..
Our aim is to prove this formula for Toeplitz operators on periodic domains, see
Theorem \ref{th8.3}, but this section is devoted to the basic definitions and
some preparatory results. 

From now on we consider Toeplitz operators $T_a$ with periodic symbols 
$a \in L^\infty (\Pi)$: we assume that 
\bea
a(z) = a(z +1) \ \ \ \mbox{for almost all} \ z \in \Pi . \label{8.8}
\eea
According to \ef{1.1a}, the Toeplitz operator $T_a$ is defined by
\bea
T_a f = P_\Pi M_a f = P_\Pi(af) , 
\label{8.10}
\eea
where $M_a $ is the pointwise multiplier. Clearly, $T_a$ is bounded on $A^2(\Pi)$,
since $a \in L^\infty(\Pi)$. 
In addition, 
we define for all $\eta \in [-\pi, \pi]$ the bounded, 
Toeplitz-type  operator $T_{a,\eta} : A_\eta^2 (\varpi) \to A_\eta^2(\varpi)$,
\bea
T_{a,\eta}  f = P_\eta (a|_\varpi f)   \label{8.12}
\eea
and the corresponding operator $\cT_a : L^2(-\pi,\pi;A_\eta^2(\varpi))
\to L^2(-\pi,\pi;A_\eta^2(\varpi))$, 
\bea
\cT_a : f(\cdot , \eta) \mapsto T_{a,\eta}  f(\cdot, \eta) ,
 \label{8.11r}
\eea
and  
the operator $\cM_a : L^2(-\pi,\pi; A_\eta^2(\varpi)) \to  L^2(-\pi,\pi; A_\eta^2(\varpi)) $
by 
\bea
\cM_a : f(\cdot , \eta) \mapsto a|_\varpi f(\cdot, \eta),
\ \ \ \ \mbox{hence, \ \  $\cT_a = \cP \cM_a$,}  \label{8.11M}
\eea
see the definition of $\cP$ in Remark \ref{rem2.3}.

We denote the spectrum of $T_{a,\eta}$ in the space $A_\eta^2 (\varpi)$ 
by $\sigma(T_{a,\eta}) $. 
The following is an immediate consequence of the definitions. 

\BEL \label{lem8.0a}
We have 
\beas
T_a f =  \sF^{-1} \cT_a  \sF f   
\eeas 
for all $f \in A^2(\Pi)$. 
\ENL

Proof. 
By the definition of the Flo\-quet-transform and \ef{8.8} we have $\sF M_a = \cM_a \sF$,
hence, 
$T_af = P_\Pi M_a f = \sF^{-1} \cP \sF M_af = 
\sF^{-1} \cP \cM_a \sF f = \sF^{-1} \cT_a \sF f.  \ \ \Box$ 

\bigskip

Given $\eta, \mu \in [-\pi,\pi]$ we denote by $J_{\eta,\mu} :
A^2(\varpi) \to A^2(\varpi)$ the operator
\beas
J_{\eta,\mu} f (z) = e^{i(\mu - \eta) z} f(z) . 
\eeas
The operator $J_{\eta, \mu}$ maps $A^2(\varpi)$ bijectively onto itself and also 
$A_\eta^2(\varpi)$ bijectively onto $A_\mu^2(\varpi)$, although it is not an isometry, 
since the factor $e^{i(\mu -\eta)z} $ is not unimodular.
Using the Taylor series  and taking into account
$|z| \leq M +1$ for $z \in \varpi$ yield 
$|1 - e^{i (\mu - \eta) z}| \leq C |\mu -\eta| $, where $C$ depends on $M$, but
this is not indicated.  Thus, 
\bea
\Vert I_\varpi - J_{\eta, \mu} \Vert_{L^2(\varpi) \to L^2(\varpi)} \leq C |\eta - \mu| \ \ \ \forall \, \mu, \eta
\in [-\pi, \pi] , \label{8.16a}
\eea
where $I_\varpi$ is the identity operator on $A^2(\varpi)$. We also denote
\bea
T_{a,\eta,\mu} = J_{\eta, \mu}^{-1}  T_{a,\mu} J_{\eta, \mu}
 = J_{\mu, \eta}  T_{a,\mu} J_{\eta, \mu} :
A_\eta^2(\varpi)  \to A_\eta^2(\varpi)
 \label{8.16ax}
\eea
and we will consider $T_{a,\eta,\mu}$ also as an operator 
$A^2(\varpi) \to A^2(\varpi)$. The inequality \eqref{8.16a} implies
\bea
\Vert T_{a,\mu} - T_{a,\eta, \mu} \Vert_{A^2(\varpi) \to A^2(\varpi)} \leq C |\eta- \mu| . 
\label{8.12d}
\eea

\BEL  \label{lem8.0} 
There exist constants $C, C' >0$ such that if  $\eta, \mu \in [-\pi,\pi]$, we have 
$\Vert P_\eta - P_\mu \Vert_{L^2(\varpi) \to L^2(\varpi)} \leq C |\eta- \mu|^{1/2}$ 
and consequently $\Vert T_{a,\eta} - T_{a,\mu} \Vert_{A^2(\varpi) \to A^2(\varpi)} \leq 
C' |\eta- \mu|^{1/2} $. 
\ENL

For simply connected domains this can be derived from formula (6.14), \cite{T1},
for the integral kernel of $P_\eta$.
The general case can be treated with the help of the isomorphism $J_{\eta,\mu}$,
as will be shown next.

\bigskip

Proof.  We fix $\mu \in [-\pi,\pi]$ and for all $\eta$ define the operator
$\widetilde P_\mu = J_{\eta,\mu} P_\eta J_{\mu,\eta}$, which is a non-orthogonal
projection from $L^2(\varpi) $ onto $A_\mu^2(\varpi)$. Formula \ef{8.16a}
implies
\beas
\Vert P_\eta J_{\mu,\eta} f - \wP f \Vert_\varpi
\leq C |\eta - \mu| \Vert f \Vert_\varpi
\eeas
for all  $f \in L^2(\varpi)$. This and another use of \ef{8.16a} yield for all
$f \in L^2(\varpi)$
\bea
& & \big| \big( f - \wP f \big| \wP f \big)_\varpi \big|
\roweq
\Big| \big( f - \wP f \big| \wP f \big)_\varpi -  
\big( J_{\mu,\eta} f - P_\eta J_{\mu,\eta} f \big| P_\eta J_{\mu,\eta} f \big)_\varpi \Big| 
\rowleq
\Big| \big( f - J_{\mu,\eta}f +P_\eta J_{\mu,\eta} f- \wP f \big| \wP f \big)_\varpi \Big|
\rowpl
\Big|  \big( J_{\mu,\eta} f - P_\eta J_{\mu,\eta} f \big| \wP f - P_\eta J_{\mu,\eta} 
f \big)_\varpi  \Big| 
\rowleq
\Big| \big( f - J_{\mu,\eta}f  \big| \wP f \big)_\varpi \Big|
+ \Big| \big( P_\eta J_{\mu,\eta} f- \wP f \big| \wP f \big)_\varpi \Big|
+ C |\eta - \mu| \Vert f \Vert_\varpi^2 
\rowleq
C' |\eta - \mu| \Vert f \Vert_\varpi^2 
, \label{8.s4}
\eea
where an inner product vanishes on the second line,  since $P_\eta$ is an orthogonal 
projection. Now, given $f \in L^2(\varpi)$, we write $f_A = P_\mu f$, $f^\perp =
f -f_A$, and, since both $P_\mu$ and $\wP$ project onto $A_\mu^2(\varpi)$, we obtain 
$(P_\mu - \wP) f = (P_\mu - \wP) f^\perp$ and thus
\bea
& & \big( (P_\mu - \wP) f \big| (P_\mu - \wP) f \big)_\varpi 
=  \big( (P_\mu - \wP) f^\perp \big| (P_\mu - \wP) f^\perp \big)_\varpi
\roweq
( \wP f^\perp | \wP f^\perp )_\varpi.
\label{8.s6}
\eea
We apply  \ef{8.s4} to $f^\perp$ and get
\bea
& & \big| ( \wP f^\perp | \wP f^\perp )_\varpi \big| =
\big| ( \wP f^\perp | \wP f^\perp )_\varpi  -  (\wP  f^\perp | f^\perp )_\varpi \big|
\rowleq 
C|\eta -\mu| \Vert f \Vert^2 , \label{8.s8}
\eea
since  here $ (\wP  f^\perp | f^\perp )_\varpi=0$ as  $ \wP  f^\perp \in 
A_\mu^2(\varpi)$ and $f^\perp$ is in the orthogonal complement. By combining \ef{8.s6} and \ef{8.s8} we obtain
$\Vert P_\mu - \wP \Vert_{L^2(\varpi) \to L^2(\varpi)}$ $ \leq 
C |\eta -\mu|^{1/2}$. This implies the lemma, since \ef{8.16a} yields 
\beas
& & \Vert P_\eta -\wP \Vert_{L^2(\varpi) \to L^2(\varpi)}  \leq \Vert( I_\varpi - J_{\eta, \mu }) P_\eta \Vert_{L^2(\varpi) \to L^2(\varpi)}  
\rowpl
\Vert J_{\eta, \mu } P_\eta ( I_\varpi - J_{\mu,\eta} ) \Vert_{L^2(\varpi) \to L^2(\varpi)}  \leq C |\eta -\mu|.
\hskip1cm \ \ \Box
\eeas

In the rest of the section we record some elementary facts in spectral theory,  
adapted to our considerations and needed later.

\BEL  
\label{lem8.0g}
If  $\eta, \mu \in [-\pi,\pi]$,  
the spectra of the operators $T_{a,\mu}: A_\mu^2(\varpi) \to A_\mu^2(\varpi) $ 
and $T_{a,\eta,\mu}: A_\eta^2(\varpi) \to A_\eta^2(\varpi) $  coincide. 
\ENL

Indeed, if  $S \in \cL (A_\mu^2(\varpi))$ and  $\lambda \in \bbC$, we have 
\bea
\big(J_{\eta,\mu}^{-1} S J_{\eta,\mu} - \lambda I_\varpi \big)^{-1} = 
J_{\eta,\mu}^{-1}(  S - \lambda I_\varpi  )^{-1} J_{\eta,\mu} ,
\label{8.11y}
\eea
provided one of the operators $J_{\eta,\mu}^{-1} S J_{\eta,\mu} - \lambda I_\varpi $
and $ S - \lambda I_\varpi  $ has a bounded inverse.

As well known, for a given Hilbert space $X$ and  $T \in \cL(X)$, the resolvent $R_\lambda(T)$ 
depends continuously and even analytically on $\lambda \in \varrho (T) $:
if $\mu \in \varrho(T)$ is fixed and $\lambda$ is such that 
say, 
\bea
|\mu -\lambda| \leq \frac12 \Vert R_\mu(T) \Vert ,\label{8.12f}
\eea
then 
$ 
R_\lambda (T) = \big( I_X - (\lambda - \mu ) R_\mu (T) \big)^{-1} R_\mu(T) .
$ 
Here, the inverse operator is defined by using the Neumann series,
which also yields an  operator norm  estimate 
\bea
\Vert R_\lambda (T) \Vert \leq (1 + 2 |\mu -\lambda| ) \Vert R_\mu(T) \Vert^2 ,
\label{8.12g}
\eea
if \eqref{8.12f} holds. 

We will need another variation of the same theme.

\BEL
\label{lem8.0c}
Let $X$ be a Hilbert space and assume  $S \in \cL(X)$ and  $\lambda \in \varrho 
(S) $. Then, $\lambda$ belongs to the resolvent set for all bounded operators 
in the set
\beas
W(S) := \big\{ T \in \cL(X) \, : \, \Vert T -S \Vert  < 
\Vert R_\lambda(S) \Vert^ {-1} \big\} \subset \cL(X).  
\eeas
Moreover, the mapping  $\cR: W(S) \to \cL(X)$, $\cR : T \mapsto R_\lambda(T)$ 
is locally Lipschitz-continuous with respect to the operator norm. 
\ENL

Indeed, the resolvent $R_\lambda (T)$ of $T \in W(S)$ is obtained for $\lambda \in \varrho 
(S) $ from the Neumann series
\bea
R_\lambda(T) = ( T - \lambda I)^{-1} =  \big(I + R_\lambda(S) ( T -S) \big)^{-1} R_\lambda(S) ,
\label{8.13a}
\eea
which obviously converges due to the definition of $W(S)$.
Moreover, applying formula \eqref{8.13a} to the given $T , \widetilde T \in  W(S)$  and 
$\lambda \in \varrho (T) $ yields 
\bea
& & \Vert R_\lambda(T) - R_\lambda(\widetilde T ) \Vert 
= \Big\Vert R_\lambda(T) - \Big( I + \sum_{n=1}^\infty ( R_\lambda(T)
(T - \widetilde T ))^n \Big)R_\lambda(T) \Big\Vert
\rowleq 
 \Vert R_\lambda(T) \Vert \sum_{n=1}^\infty  \Vert R_\lambda(T) \Vert^n
\Vert T- \widetilde T \Vert^n  \leq C_{T,\lambda} \Vert \widetilde T - T\Vert ,
\eea
provided $\widetilde T$ belongs to a small enough neighborhood of $T$ in the 
space $\cL(X)$.

\section{Essential spectrum of a Toeplitz operator with periodic symbol}  
\label{sec9}

Let us state the main result of this paper. Here, the periodic domain $\Pi$ is as 
in the beginning of Section \ref{sec2}, the symbol $a \in L^\infty(\Pi)$ of the Toeplitz-
operator $T_a$ is periodic according to \ef{8.8},  the operators $T_{a,\eta}$ have been   
defined in \ef{8.12}, and $\sigma(T_{a,\eta})$ denotes the spectrum of $T_{a,\eta}$
as an operator in the space $A_\eta^2(\varpi)$.

\BET  \label{th8.3}
The essential spectrum of the Toeplitz-operator $T_a:A^2(\Pi) \to A^2(\Pi)$ can be described by the
formula
\bea
\sigma_{\rm ess} ( T_a )  = \bigcup_{\eta \in [-\pi,\pi]} \sigma(T_{a,\eta}). 
\label{8.40}
\eea
Moreover, there holds $\sigma( T_a )= \sigma_{\rm ess} ( T_a )$.
\ENT

We will give the proof in several steps. Let us start with the following observation. 

\BEL \label{lem8.1}
The set 
$
\Sigma := \bigcup_{\eta \in [-\pi,  \pi]} \sigma(T_{a,\eta})  
$ 
is closed. 
\ENL

Proof. We first fix the constants 
\beas
& & D_0  = 1/(8D_1^2  D_2) > 0 , \ \ \ \mbox{where} \ D_1 = \sup_{\eta,\mu \in [-\pi,\pi]} \Vert J_{\eta,\mu} \Vert_{A^2(\varpi)\to A^2(\varpi)} ,
\eeas
and $D_2$ is the largest of the constants $C, C'$ appearing in \ef{8.16a}, \ef{8.12d} and Lemma \ref{lem8.0}.

If $\Sigma$ were not closed, we would find  $\lambda \in 
{\rm cl}\,(\Sigma) \setminus \Sigma$ and sequences $(\eta_k)_{k=1}^\infty 
\subset [-\pi, \pi]$ and $(\lambda_k)_{k=1}^\infty$ such that $\lambda_k \in 
\sigma (T_{a,\eta_k} )$ and $\lambda_k \to \lambda$ as $k \to \infty$. Due to the
compactness of the  $\eta$-interval, by passing to a subsequence we may assume 
that for some $\eta$, fixed from now on, we have  $\eta_k \to \eta$ as $k \to  \infty$. 
Since $\sigma(T_{a,\eta})$ is closed, we find $0< \delta <1$ such that 
$ 
{\rm dist}\,( \lambda, \sigma(T_{a,\eta}) ) \geq \delta . 
$ 
We now fix $k$ large enough such that 
$ 
|\lambda_k - \lambda| \leq \delta /2$, which implies  that 
$\lambda_k \in \varrho(T_{a,\eta})$, 
and such that 
\bea
|\eta -\eta_k |^{1/2} \leq \min \big( 1, D_0 \Vert R_\lambda (T_{a,\eta} )
\Vert_{A_\eta^2(\varpi)  \to A_\eta^2(\varpi)}^{-2} \big) ,
\label{8.20p}
\eea
and consider $S:= T_{a,\eta_k,\eta}: A_{\eta_k}^2(\varpi) \to A_{\eta_k}^2(\varpi)$. 
According to Lemma \ref{lem8.0g}, the spectra of $T_{a,\eta}$  in 
$A_\eta^2(\varpi)$ and $S$ in $A_{\eta_k}^2(\varpi)$ coincide, 
hence we have $\lambda_k \in \varrho(S)$. Moreover, by \eqref{8.12d}
and Lemma  \ref{lem8.0},
\bea
& & \Vert S- T_{a,\eta_k} \Vert_{A^2(\varpi) \to A^2(\varpi)} 
\rowleq 
\Vert S- T_{a,\eta} \Vert_{A^2(\varpi) \to A^2(\varpi)} 
+ \Vert T_{a,\eta}- T_{a,\eta_k} \Vert_{A^2(\varpi) \to A^2(\varpi)}
\leq 2 D_2 |\eta -\eta_k|^{1/2} .  \label{8.20q}
\eea
Finally, we have, by \ef{8.11y} and \eqref{8.12g} 
\beas
&&\Vert R_{\lambda_k}(S) \Vert_{A_{\eta_k}^2(\varpi)\to A_{\eta_k}^2(\varpi)} 
= \Vert J_{\eta_k,\eta}^{-1} R_{\lambda_k}(T_{a,\eta} )J_{\eta_k,\eta} 
\Vert_{A_{\eta_k}^2(\varpi)\to A_{\eta_k}^2(\varpi)}
\rowleq 
\Vert J_{\eta_k,\eta}^{-1} \Vert_{A_\eta^2(\varpi)\to A_{\eta_k}^2(\varpi)} \, \Vert J_{\eta_k,\eta} \Vert_{A_{\eta_k}^2(\varpi)\to A_\eta^2(\varpi)} 
\nonumber  \\
& & \ \ \times
(1 + 2 |\lambda_k - \lambda|) 
\Vert R_{\lambda}(T_{a,\eta} ) \Vert_{A_\eta^2(\varpi)\to A_\eta^2(\varpi)}^2
\rowleq 
\Vert J_{\eta_k,\eta}^{-1} \Vert_{A^2(\varpi)\to A^2(\varpi)} \, \Vert J_{\eta_k,\eta} \Vert_{A^2(\varpi)\to A^2(\varpi)}
\nonumber  \\
& & \ \ \times
(1 + 2 |\lambda_k - \lambda|) 
\Vert R_{\lambda}(T_{a,\eta} ) \Vert_{A_\eta^2(\varpi)\to A_\eta^2(\varpi)}^2
\rowleq 
2 D_1^2   \Vert R_{\lambda}(T_{a,\eta} ) 
\Vert_{A_\eta^2(\varpi)\to A_\eta^2(\varpi)}^2 . 
\eeas
Hence, by \eqref{8.20q}, \ef{8.20p} and the choice of the various constants
in the beginning of the proof, 
\beas
& & \Vert S- T_{a,\eta_k} \Vert_{A_{\eta_k}^2(\varpi)\to A_{\eta_k}^2(\varpi)}
\leq
\Vert S- T_{a,\eta_k} \Vert_{A^2(\varpi)\to A^2(\varpi)}
\rowleq 
\frac12 \Vert R_{\lambda_k}(S) 
\Vert_{A_{\eta_k}^2(\varpi)\to A_{\eta_k}^2(\varpi)}^{-1}.
\eeas
Consequently, Lemma \ref{lem8.0c} applies and shows that 
$\lambda_k \in \varrho(T_{a,\eta_k})$, which is a contradiction. 
\ \ $\Box$

 \BEC
\label{th8.2}
The set $\Sigma $ is contains the spectrum $\sigma (T_a)$ and thus also
the essential spectrum $\sigma_{\rm ess} (T_a)$  of the operator $T_a$. 
\ENC

Proof. Let us fix a number $\lambda \in \bbC$ with $\lambda \notin \Sigma$. 
Then, for every $\eta \in [-\pi,\pi]$, there exists a bounded inverse
$R_{\eta, \lambda} : A_\eta^2(\varpi) \to A_\eta^2(\varpi)$ of the operator 
$T_{a,\eta} - \lambda I_\varpi$. Let us show that the operator norm of 
$R_{\eta, \lambda}$ in the space $A_\eta^2(\varpi)$  depends continuously on $\eta$.
Fixing $\mu \in [-\pi,\pi]$, we note by \ef{8.11y}, \ef{8.16a} that
\bea
& & \Vert ( T_{a,\eta,\mu }- \lambda I_\varpi)^{-1} \Vert_{A_\eta^2(\varpi) \to A_\eta^2(\varpi)} 
\roweq 
\Vert J_{\mu,\eta} (T_{a,\mu} - \lambda I_\varpi )^{-1} J_{\eta,\mu} 
\Vert_{A_\eta^2(\varpi) \to A_\eta^2(\varpi)} \to
\Vert (T_{a,\mu} - \lambda I_\varpi )^{-1} \Vert_{A_\mu^2(\varpi) \to A_\mu^2(\varpi)}  \label{spec}
\eea
as $\eta \to \mu$ in $[-\pi,\pi]$. 
But \ef{8.12d} and Lemma \ref{lem8.0} also yield
\beas
\Vert T_{a,\eta,\mu } - T_{a,\eta} \Vert_{A_\eta^2(\varpi) \to A_\eta^2(\varpi)}
\to 0
\eeas
as $\eta \to \mu$, hence, by Lemma \ref{lem8.0c},
\beas
& & \Vert ( T_{a,\eta }- \lambda I_\varpi)^{-1} \Vert_{A_\eta^2(\varpi) \to A_\eta^2(\varpi)} 
\to
\Vert (T_{a,\eta, \mu} - \lambda I_\varpi )^{-1} \Vert_{A_\eta^2(\varpi) \to A_\eta^2(\varpi)}  
\eeas
as $\eta \to \mu$ in $[-\pi,\pi]$. This together with \ef{spec}
prove the claim. 

The compactness of the $\eta$-interval implies that the operator norm
$R_{\eta, \lambda} : A_\eta^2(\varpi) \to A_\eta^2(\varpi)$ has a uniform
upper bound for all $\eta \in [-\pi,\pi]$. 
%
%
Consequently, the operator $\cR_\lambda : f(\cdot ,\eta) \mapsto R_{\eta, \lambda} 
f(\cdot,\eta)$ is bounded in $L^2(-\pi,\pi ; A^2_\eta(\varpi))$. Denoting
by $\cI$ the indentity operator on $L^2(-\pi,\pi ; A^2_\eta(\varpi))$, it is 
clear (see \ef{8.11r}) that $\cR_\lambda$ is the inverse operator of 
$\cT_a - \lambda \cI$ in the space  $L^2(-\pi,\pi ; A^2_\eta(\varpi))$ . 
Hence, 
$\sF^{-1} \cR_\lambda \sF$ is a bounded inverse of $T_a - \lambda 
I_\Pi = \sF^{-1} ( \cT_a - \lambda\cI) \sF$, and $\lambda$ belongs to the 
resolvent set of $T_a$. 
\ \ $\Box$

\bigskip

The proof of the converse relation needs the following simple observation.

\BEL \label{lem8.3}
For every $n\in\bbN$, we denote $\varphi_n(z)= e^{-n^{-2}z^2}$, $z \in \Pi$. 
Then, there holds $|\varphi_n(z)| \leq C$ for all $z \in \Pi$, $n \in \bbN$,
and also 
\beas
\varphi_n f \to f \ \ \mbox{in $ A^2(\Pi) $  as $n \to \infty$,}  
\eeas
for every $f \in A^2(\Pi)$.
\ENL

Proof. The  uniform boundedness of $\varphi_n$ follows from the boundedness of the 
domain $\Pi$ in the imaginary direction.
Moreover, if $f \in A^2(\Pi)$ and  $\varepsilon> 0$ are given, we find $M > 0$ such that  $\int_{ \Pi \cap \{ |{\rm Re}z| \geq M  \} }
|f|^2 dA \leq \varepsilon$.
Then, we pick up $n \in \bbN$ so large that 
$ 
| 1-   \exp ( - n^{-2} z^2 ) | \leq
\varepsilon
$ 
for all $z$ with $ |{\rm Re}z| \leq  M$ and get
\beas
\Vert \varphi_n f - f   \Vert_\Pi^2 \leq 
\!\!\!\!\!\! \!\!\!\!\!\!
\int\limits_{ \Pi \cap \{ |{\rm Re}z| \geq M  \} }
\!\!\!\!\!\! \!\!\!\!\!\!
C |f|^2 dA 
+ \!\!\!\!\!\! \!\!\!\!\!\!
\int\limits_{ \Pi \cap \{ |{\rm Re}z| \leq M  \} }
\!\!\!\!\!\! \!\!\!\!\!\!
\big| 1-   e^{- n^{-2} z^2 } \big| |f |^2 dA 
\leq C\varepsilon +\varepsilon \Vert f \Vert_\Pi^2 .  \ \ \Box
\eeas

\bigskip

Proof of Theorem \ref{th8.3}. In view of Corollary \ref{th8.2}, it suffices to show 
that every $\lambda$ which belongs to $\sigma(T_{a,\mu})$ for some 
$\mu \in [-\pi, \pi]$ is a point in the essential spectrum of $T_a$. Thus, let us fix such $
\lambda$ and $\mu$. We aim to use the Weyl criterion of Lemma \ref{lemW} in order
to show that $\lambda \in \sigma_{\rm ess} (T_a)$. 

Let $\varepsilon > 0$ be arbitrary. Now, $\lambda$ is
an eigenvalue or a point in the essential spectrum of $T_{a,\mu}$ so that for
example by Lemma \ref{lemW} we can 
find a near eigenfunction $g \in A_{\mu}^2(\varpi)$ with 
$\Vert g \Vert_\varpi = 1$ such that
\bea
\Vert T_{a,\mu} g - \lambda g \Vert_\varpi \leq \varepsilon. \label{8.42}
\eea
By a suitable approximation, we may assume that in fact $g$ has an analytic
extension in a neighborhood of $\varpi$ and satisfies the quasiperiodic
boundary condition \ef{1.15} (see the definition of $A_{\eta, {\rm ext}}^2(\varpi)$ in 
Section \ref{sec2}). 

We set for every $n \in \bbN$
\bea
G_n = \sF^{-1} (\cX_{\mu,n} g ) ,   \label{8.43}
\eea
where  $  (\cX_{\mu,n}g) (z,\eta) := g \otimes\cX_{\mu,n}  (z,\eta) =  
g(z)\cX_{\mu,n}(\eta) \in L^2(-\pi,\pi;A_\eta^2(\varpi))$ with 
\begin{equation}
\cX_{\mu,n}(\eta) = \left\{
\begin{array}{ll}
\sqrt{n} , \ \ \ & |\mu - \eta| \leq \frac1{2n} \\
0, &\mbox{for other} \ \eta \in [-\pi,\pi]. 
\end{array}
\right. \label{8.43a}
\end{equation}
Since $ \Vert  g \Vert_\varpi=1$, we have 
$
\frac12 \leq \Vert \cX_{\mu,n} g \Vert_{L^2 (-\pi,\pi; L^2(\varpi)) } \leq 1,
$
and thus the isometric property of the Floquet transform implies
\bea
G_n = \sF^{-1} ( \cX_{\mu,n} g )  \in  A^2(\Pi) \ \ \ \mbox{and} \ \ \ 
\frac12 \leq \Vert G_n  \Vert_\Pi \leq 1.  \label{8.43p}
\eea
Next, we recall  Lemma \ref{lem8.3}
and choose for every $n$ the number $\ell(n) \in \bbN$ such that $\ell(n+1) > \ell(n)$
and such that $\Vert G_n -  \varphi_{\ell(n)}  G_n  \Vert_\Pi \leq 1/(n+2)$. By 
reindexing  $\varphi_{\ell(n)} \to \varphi_n$ we thus have, for all $n \in \bbN$,   
\bea
\Vert G_n -  \varphi_{n}  G_n  \Vert_\Pi \leq \frac{1}{n +2}
\label{8.43pp}
\eea
and then define 
\bea
f_n = \varphi_{n} G_n \in A^2(\Pi) \ \ \
\Rightarrow \ \ \ \frac14 \leq \Vert f_n \Vert_\Pi \leq 2, \label{8.43ppp}
\eea
see \ef{8.43p}.  

We will obtain the desired Weyl singular sequence by taking translations of these
functions, but let us first show the crucial estimate 
\bea
\Vert T_a f_n - \lambda f_n \Vert_\Pi \leq C \varepsilon. \label{8.39k} 
\eea 
To see this, let $n  \geq 1/\varepsilon^2$ in the following.
The boundedness of the operator $T_a$ in $L^2(\Pi)$ and 
\ef{8.43pp}, \ef{8.43ppp} yield
\bea
& & \big\Vert T_a ( f_n -  G_n) \big\Vert_\Pi 
< C \varepsilon. \label{8.56a}
\eea
Moreover, 
\bea
T_a G_ n &= &  \sF^{-1} \cP \cM_a (\cX_{\mu,n} g  ) 
= \sF^{-1}   \big(  \cX_{\mu,n }(\eta) P_\eta ( a g  )(z)  \big) 
\roweq
\sF^{-1}  \Big(   \cX_{\mu,n } (\eta) \big( P_\eta ( a g  ) (z)
- P_\mu (ag)(z) \big) \Big)  
\rowpl
\sF^{-1}   \big( \cX_{\mu,n } (\eta) P_\mu (ag)(z)  - \lambda  \cX_{\mu,n }(\eta) g(z) \big)  
+ \sF^{-1}   \big( \lambda  \cX_{\mu,n } g \big)  
\nonumber \\
&=:& 
\Psi_1 + \Psi_2  +\lambda G_n .  \label{8.56}
\eea
To estimate $\Psi_1$, 
Lemma \ref{lem8.0}  states that  $\Vert P_\eta - P_\mu \Vert_{A^2(\varpi) \to 
A^2 (\varpi)} \leq C|\eta- \mu|^{1/2}$, hence, by \ef{8.43a} and the boundedness
of the symbol $a$, 
\bea
& & \big\Vert    \cX_{\mu,n } \big( P_\eta ( a g  ) - 
P_\mu (ag)  \big) \big\Vert_{L^2(-\pi,\pi;L^2(\varpi))}^2
\leq 
C_a \int\limits_{\mu - 1/(2n)}^{\mu+1/(2n)} n|\eta -\mu| \Vert g \Vert_\varpi^2 d\eta
\leq 
\frac{C'_a}{n} .  \label{8.56w}
\eea  
The near eigenfunction property \ef{8.42} yields a similar estimate for $\Psi_2$:
\beas
& & \big\Vert    \cX_{\mu,n } \big( P_\mu ( a g  ) - 
\lambda g  \big)\big\Vert_{L^2(-\pi,\pi;L^2(\varpi))}^2
\leq 
C \int\limits_{\mu - 1/(2n)}^{\mu+1/(2n)} n \varepsilon^2  d\eta
\leq C' \varepsilon^2. 
\eeas
From this and \ef{8.56w} we obtain $\Vert \Psi_1 \Vert_\Pi + 
\Vert \Psi_2 \Vert_\Pi \leq C \varepsilon$. Taking into account \ef{8.43pp}, \ef{8.43ppp}, 
\ef{8.56a},  and \ef{8.56} gives \ef{8.39k}. 

We finally prove  that the translated functions $h_n = f_n \circ \sft_n
\Vert f_n\Vert_\Pi^{-1}$,
where $\sft_n(z) = z -  m(n) $, form a Weyl singular sequence as defined around \ef{8.39},
if the integers $m(n)$ are chosen suitably. This will be done soon, but we first remark
that since the symbol $a$ is 1-periodic with respect to the real coordinate and the Bergman
projection $P_\Pi$ commutes with the translations $\tau_n$, see Remark \ref{rem2.3},
there holds the norm invariance
\beas
\Vert T_a f_n - \lambda f_ n \Vert_\Pi = 
\Vert (T_a f_n ) \circ \sft_n - \lambda  f_n \circ \sft_n \Vert_\Pi
= \Vert T_a h_n - \lambda h_n \Vert_\Pi  \, \Vert f_ n \Vert_\Pi 
\eeas
and thus  \ef{8.39k} implies the convergence \ef{8.39}, in view of \ef{8.43ppp}.

There remains to show that the sequence $(h_n)_{n=1}^\infty$ does not have
convergent subsequences for a proper choice of the numbers $m(n)$ in the 
definition of the translations $\sft_n$. 
For every $n$, let the integers $M(n) > M(n-1) > \ldots  > 0$ be such that 
$|\varphi_{n}(z) | \leq 1/n$ for all $z$ with $|{\rm Re}\,z| \geq M(n)$ and also
such that
\bea
\int\limits_{\Pi \cap \{ |{\rm Re}z| \leq  M(n)  \} }
\!\!\!\!\!\! \!\!\!\!\!\! | \varphi_{n} G_n|^2 dA \geq 
\Vert\varphi_{n}  G_n  \Vert_\Pi^2 - \frac{1}{n} \geq \Vert G_n \Vert_\Pi^2 - \frac2n
\geq \frac14 - \frac2n,  \label{8.60}
\eea
see \ef{8.43p}, \ef{8.43pp}, and then choose
by induction, $m(1) = M(1)$ and $m(n) = m(n-1) + 2 M(n) $. Next, we denote 
$$
\Pi_n = \big\{ z \, : \, {\rm Re} \, z \in \Gamma_n:=  [m(n)- M(n) , m(n) + M(n) 
\big\}
$$ 
for all $n$, and observe that for the indices $n \not= \ell$  there always holds
$\Pi_n \cap \Pi_\ell = \emptyset$. This is so, since, e.g. in the case $n > \ell$,
the lower bound of the interval $\Gamma_n$ equals
$m(n) -M(n) = m(n-1) + M(n) > m(\ell) + M(\ell)$ which is the upper bound of the 
interval $\Gamma_\ell$. Thus, for all $n \not= \ell$ we obtain by the definition
of $\sft_\ell$ and the choice of the number $M(\ell)$ that 
\beas
|\varphi_\ell \circ \sft_\ell  (z) | \leq \frac1\ell  \ \ \ \forall \, z .
\in \Pi_n
\eeas
Finally, \ef{8.43p}, \ef{8.43pp}, \ef{8.43ppp}, \ef{8.60}, imply for 
large enough $n$
\beas
& & \Vert h_n - h_\ell \Vert_\Pi^2  \geq \int\limits_{\Pi_n} |h_n -h_\ell|^2 dA
\rowgeq 
\int\limits_{\Pi_n} \frac{|\varphi_n \circ \sft_n \,  G_n \circ \sft_n|^2}{
\Vert f_n \Vert^2}  dA -
\int\limits_{\Pi_n} \frac{|\varphi_\ell \circ \sft_\ell \,  G_\ell \circ \sft_\ell|^2}{
\Vert f_\ell \Vert^2} dA    
\rowgeq
\!\!\!\!\!\!
\int\limits_{\Pi \cap \{ |{\rm Re}z| \leq  M(n)  \} }
\!\!\!\!\!\! \!\!\!\!\!\!  \frac{|\varphi_{n} G_n|^2 }{ \Vert f_n \Vert^2} dA
- \frac1\ell \int\limits_{\Pi} \frac{ | G_\ell |}{ \Vert f_\ell \Vert^2}^2 dA    
\geq
\frac1{16} - \frac{1}{2n}  - \frac{16}{\ell}.
\eeas
Thus, no convergent subsequences exist. 
This completes the proof of the theorem, since we have shown that $\lambda 
\in \sigma_{\rm ess} (T_a)$. 
 \ \ $\Box$

\section{Toeplitz operator with prescribed essential spectrum:
statement and definitions.} \label{sec9R}

We present  examples of Toeplitz operators with interesting essential 
spectra.

\BET  \label{cor9.2}
Given $K \in \bbN$ and any finite sequence of distinct real numbers $x_1 > \ldots >  x_K$
one can find $\delta > 0$ such that for all small enough $\varepsilon > 0$, there exists a bounded 
Toeplitz-operator $T_a:A^2(\bbD) \to A^2(\bbD)$ with a real valued symbol
$a \in L^\infty(\bbD)$ and the properties 
\bea
& & \sigma_{\rm ess} (T_a)  \cap  B(x_n, \varepsilon)
\not= \emptyset 
\ \forall \, n \leq K  \label{9.2bb}
\eea
and
\bea
& & {\rm dist} \, \big( \sigma_{\rm ess} (T_a) \cap U_\varepsilon ,
\sigma_{\rm ess}(T_a) \setminus U_\varepsilon\big) \geq \delta, 
\ \ \mbox{where} \ U_\varepsilon = \bigcup_{n \leq K} B(x_n, \varepsilon).
\label{9.2bbb}
\eea
\ENT

This result is a consequence of the  following theorem on the essential spectra of
Toeplitz operators on periodic domains. The idea of the proof of Theorem \ref{th9.1}
is discussed below Definition \ref{def5.3}.

\BET
\label{th9.1}
Let $b \in L^\infty(\bbD)$ be a real valued symbol such that its support is 
contained in $\bbD_R$ for some $R< 1$. Denote by $\Sigma =  \{ 0 \} \cup 
\bigcup_{n=1}^\infty \{\lambda_n \}$ the spectrum of $T_b$, where the eigenvalues 
$\lambda_n$ are indexed such that
\bea
|\lambda_1| \geq | \lambda_2| \geq |\lambda_3| \geq \cdots 
\label{9.1}
\eea 
and such that any $\lambda_n$ appears as many times as its  multiplicity is. 

If $N \in \bbN$ and $\delta > 0$  are arbitrary numbers  such that $|\lambda_N| > |
\lambda_{N+1}| + 2 \delta$, then for every small enough $\varepsilon>0$ 
there exists a simply connected periodic domain $\Pi \subset \bbC$ as in \ef{1.2}
and a bounded Toeplitz operator $T_{\mathsf a} :  A^2(\Pi) \to A^2(\Pi)$ 
such that 
\bea  
& & \sigma_{\rm ess} (T_{\sfa})  \cap  B(\lambda_n, \varepsilon)
\not= \emptyset 
\ \forall \, n \leq N , \label{9.2b}
\eea
and
\bea
& & {\rm dist} \, \big( \sigma_{\rm ess} (T_\sfa) \cap G_\varepsilon ,
\sigma_{\rm ess} (T_\sfa)\setminus G_\varepsilon\big) \geq \delta, 
\ \ \mbox{where} \ G_\varepsilon = \bigcup_{n \leq N} B(\lambda_n, \varepsilon).
\label{9.2a}
\eea
\ENT

{\it Remarks.} 
$1^\circ$. The assumptions on the symbol $b$ imply that the 
operator $T_b$ is compact (e.g. \cite{Zh}, Theorem 7.8) and self-adjoint in $A^2(\bbD)$ 
(easy to check directly) and consequently, its spectrum is as presented above. 
The connection between the symbols $b$ and $\sfa$ is very simple: $\sfa$ is a 
periodic symbol defined by infinitely many "translated copies" of $b$, see \ef{9.11a} for 
the exact definition.

$2^\circ$. If $N$ is such that $|\lambda_N| > |\lambda_{N+1}|$, then \ef{9.2b}
and \ef{9.2a} imply that, for a sufficiently
small $\varepsilon$, the essential spectrum $\sigma_{\rm ess} (T_\sfa)$ contains at 
least as many components as there are distinct points in the set
$\{ \lambda_1, \ldots \lambda_N\}$. The components are close to these points, according
to \ef{9.2b}. We actually expect that these components (and the components of the 
essential spectrum of the operator $T_a$ of 
Theorem \ref{cor9.2}) consist of continua instead of
some number of separate points, although we do not have a proof for this assertion. The corresponding question
in the elliptic pde-theory is actually quite deep and open in the case of general
second order periodic elliptic operators, see Section 6.3 and in particular Conjecture 6.13 in \cite{Ku}, also \cite{KT}.  

$3^\circ$. The support of the symbol $a \in L^\infty(\bbD)$ will be disconnected in the 
following  construction, but is is clear that adding to $a$ for example a symbol
$ e^{ - 1/ (1 - |z|) }$ 
induces 
only a compact perturbation of the Toeplitz operator and thus does not change
the relations \ef{9.2bb} or \ef{9.2bbb}. And, the support of the modified symbol
would be connected. 

\bigskip

Proof of Theorem \ref{cor9.2}, assuming Theorem \ref{th9.1} holds. Let $K$ and the numbers
$x_1 , \ldots , x_K \in \bbR $ be given. Recall that if $b : \bbD \to \bbR$ is a bounded, 
radial symbol, then the Toeplitz-operator $T_b: A^2(\bbD) \to A^2 (\bbD)$ equals the
Taylor coefficient multiplier $T_b : \sum_{n=0}^\infty f_nz^n \mapsto
\sum_{n=0}^\infty b_n f_n z^n$, where  $b_n = \pi^{-1}(n+1)\int_0^1 b(r)
r^{2n+1} dr $, see 
\ef{Taylor}. Since the monomials $r^{2n+1}$, $n=1, \dots, K$ of a real variable $r 
\in [0,\frac12]$ form a linearly independent set in the Hilbert-space 
$L^2(0,\frac12)$, one can find a function $b \in L^2(0,\frac12)$ such that 
\bea
\int\limits_{0}^{1/2} r^{2n+1} b(r) dr = \frac\pi{n+1} x_n   \label{9.0a}
\eea 
for all $n=1, \dots, K$. We may assume that $b$ is a bounded function, and we extend $b$ 
as null onto the interval $(\frac12 , 1)$ and then radially to the open  unit disc. 
Then, 
we fix $K_0 \in \bbN$, $K_0 \geq K$ such that 
\bea
\frac{n+1}\pi  \int\limits_{0}^{1/2} r^{2n+1} b(r) dr  \leq \frac\mu{2},
\ \ \ \mbox{where} \ \mu =  \min\{ |x_j| \, : {j=1, \ldots, K} \}, \label{9.0b}
\eea
for all $n > K_0$. This is possible since $\int_0^{1/2} r^{2n+1} b(r) dr \leq 
C 2^{-n}$. 

Since the monomials $z^n$, multiplied by suitable normalization constants, form an 
orthonormal basis of $A^2(\bbD)$, it is now plain that the operator $T_b:\sum_{n=0}^\infty 
f_nz^n  \mapsto \sum_{n=0}^\infty b_n f_n z^n$ 

\smallskip

\noindent $(i)$ has eigenvalues $x_1, \dots , x_K$, which form the set $X \subset \bbR$
(see \ef{9.0a} and the choice of the numbers $b_n$), and,

\smallskip

\noindent $(ii)$ in the case  $K_0 > K$ also has a finite set  of eigenvalues 
\beas
y_n := \frac{n+1}\pi \int\limits_{0}^{1/2} r^{2n+1} b(r) dr, \ \ n= K+1, \ldots, K_0,
\eeas
some of which may coincide with the numbers $x_n$. We denote
\beas
Y = \{ y_n \, : \, n= K+1, \ldots, K_0, \ y_n \notin X \}   
\eeas
and $\varrho = {\rm dist} (Y , X ) > 0$. If $K_0 = K$, we put $Y = \emptyset$.
Also, 

\smallskip

\noindent $(iii)$ the rest of the spectrum of $T_b$ consists of 0 and eigenvalues $\lambda_n$ which are contained in the interval $[ -\mu/2 , \mu/2]$. This is so since
$T_b$ is obviously a compact operator in $A^2(\bbD)$ and \ef{9.0b} holds.

\smallskip

We now apply Theorem \ref{th9.1} for this $b$. Denoting the eigenvalues of $T_b$ 
(as in \ef{9.1}) by $\lambda_n$ , $n\in \bbN$, we first pick up $N \in \bbN$ such that $N \geq K$ and
$\mu/4 \geq |\lambda_N| > |\lambda_{N+1} | + 2 \delta$ for some $\delta>0$.
By diminishing $\delta$, if necessary, we assume  $4\delta < \varrho= {\rm dist} (Y , X )$.
Note that all eigenvalues  $\lambda_n$, which belong to the set $X$ (see $(i)$), must 
have index $n < N$, by the choice of $N$ and $\mu$. 
Let $\varepsilon > 0$ be given as in the assumption of Theorem \ref{cor9.2}: 
if necessary, we diminish it so as to satisfy $\varepsilon <  \delta $.
Summing up, there holds 
\bea
\varepsilon <  \delta  < \min \Big( \frac\varrho4 , \frac\mu8, \frac12( |\lambda_N| - |\lambda_{N+1} 
\Big).  \label{U2}
\eea
Then, we fix the simply connected periodic domain $\Pi$ and the symbol  $\sfa \in 
L^\infty(\Pi)$  as given in Theorem \ref{th9.1} for this $b$ and $\varepsilon$.

We apply Lemma \ref{lem1.5} with the Riemann mappings $\psi: \bbD \to \Pi$, $\varphi = 
\psi^{-1}$, and define 
the Toeplitz-operator $T_a$ on $A^2(\bbD)$ with $a = \sfa \circ \psi$ and with the same 
spectral properties  \ef{9.2b}, \ef{9.2a} as the operator $T_\sfa$ on $A^2(\Pi)$. 
Now, relation \ef{9.2bb}  follows from \ef{9.2b}, since all numbers $x_n$,
$n=1,\ldots,N$, are contained in the set of the eigenvalues $\lambda_n$ 
in \ef{9.2b}. Moreover, \ef{9.2bbb} is obtained by observing that if
$\lambda \in \sigma_{\rm ess} (T_\sfa)$  does not belong to the 
set $U_\varepsilon = \cup_{n=1}^K B(x_n, \varepsilon)$,  then, due to \ef{9.2b} and \ef{9.2a}, it 
must either $(I)$ belong   to $B(\lambda_n,\varepsilon)$, where $\lambda_n \in Y $, or
$(II)$ belong to $B(\lambda_n,\varepsilon)$, where $n \leq N$ but $\lambda_n \notin X \cup Y $, or $(III)$ it does not belong to the set $G_\varepsilon$ given in \ef{9.2a}. In the  case $(I)$
we have by $(ii)$ dist\,$(\lambda_n, x_\ell) \geq \varrho $ for all $\ell = 1, \ldots,K$, hence, by \ef{U2}, 
\beas
{\rm dist}\, \big( \lambda, B(x_\ell, \varepsilon ) \big) \geq {\rm dist}\, 
\big( B(\lambda_n,\varepsilon) , B(x_\ell, \varepsilon \big) 
\geq \varrho -  2\varepsilon \geq \delta  \  \forall \, \ell = 1, \ldots,K ,
\eeas
hence, dist$\, ( \lambda, U_\varepsilon ) \geq \delta $.
In the  case $(II)$ we have $|\lambda_n| \leq \mu / 2$ by $(iii)$, thus
dist\,$(\lambda_n, x_\ell) \geq \mu /2  $ for all $\ell = 1, \ldots,K$,
and thus obtain in the same way 
$ 
{\rm dist}\, \big( \lambda, U_\varepsilon )  \geq \delta  
$ 
from \ef{U2}. 
In the case $(III)$ we get directly from  \ef{9.2bbb} that dist\,$( \lambda, U_\varepsilon ) \geq {\rm dist} \, ( \lambda, G_\varepsilon ) \geq \delta$, since
$U_\varepsilon \subset G_\varepsilon$. 
We thus see that in all cases the distance of $\lambda$ from $U_\varepsilon$ is 
at least $\delta$, hence \ef{9.2bbb} holds.   \ \ $\Box$

\bigskip

We proceed with the proof of Theorem \ref{th9.1}. Let us first 
present the main definitions, see the attached figure. 

\begin{tikzpicture}

\filldraw[fill=lightgray,ultra thin] (-4.65,-0.1) rectangle (7.65,0.1);

\draw [->,thin] (-5,0) -- (9,0);

\node (A) at (-2.7,-2)  {$- 1$} ;
\node (B) at  (0.3,-2) {$0$} ;  
\node (C) at  (3.3,-2)  {$1$} ;
\node (D) at  (6.3, -2)  {$ 2$} ;
\node (E) at  (8.5,-0.4)  {$\ \ {\rm Re}$} ; 
      
\draw [->,thin] (0,-3.3)  -- 
     (0,-2.7)  node[anchor=east] {$-1$} -- 
       (0,2.7)  node[anchor=east] {$1$} --
     (0,3.2)  node[anchor=west] {${\rm Im}$}-- (0,3.7);

\draw [dashed,thin] (1.5,-3.3)  -- (1.5, -2)  
 node[anchor=west] {$1/2$}  -- (1.5,3.3)  ;
 
\draw [dashed,thin] (-1.5,-3.3)  -- (-1.5, -2)  
 node[anchor=west] {- $1/2$}  -- (-1.5,3.3)  ;
     
\draw[thick] (-4.5,0.1)  --  (7.5,0.1);  
\draw[thick] (-4.5,-0.1)  --  (7.5,-0.1);  

\filldraw[fill=lightgray,thick](-3,0) circle (1.4cm) ; 
\filldraw[fill=lightgray,thick] (0,0) circle (1.4cm) ; 
\filldraw[fill=lightgray,thick] (3,0) circle (1.4cm) ; 
\filldraw[fill=lightgray,thick] (6, 0) circle (1.4cm) ;

\filldraw[fill=gray] (0,0) circle (1cm);

\fill[fill=lightgray] (4.3,-0.1) rectangle (4.8,0.1) ;  
\fill[fill=lightgray] (7.3,-0.1) rectangle (7.8,0.1) ;  
\fill[fill=lightgray] (-4.8,-0.1) rectangle (-4.3,0.1) ;  
\fill[fill=lightgray] (1.5,-0.1) rectangle (1.8,0.1) ;  
\fill[fill=lightgray] (-1.8,-0.1) rectangle (-1.5,0.1) ;  

\draw[thick] (-1.5,-0.1) -- (-1.5,0.1);

\draw[thick] (1.5,-0.1) -- (1.5,0.1);

\draw[dashed, very thin] (4.6,-0.1) -- ( 5.3,-0.1);
\draw[dashed, very thin] (4.6,0.1) -- ( 5.3,0.1); 
\node (F) at (5.4,0)  {$2h$} ;

\draw[step=3cm,gray, very thin] (-4.3,-3.3) grid (8.3,3.3);


\draw[densely dotted,thick] (0.7,0.95) -- (2,2) ; 
\draw[densely dotted,thick] (0.7,0.4) -- (2,2)  
node[anchor=west] {$\varpi_0$};

\draw[densely dotted,thick] (0,1.4) -- (1,3.5)  
node[anchor=west] {$R_0$};

\draw[densely dotted,thick] (0,1) -- (1.5,2.7)  
node[anchor=west] {$R R_0$}; 

\draw[densely dotted,thick] (-0.5,0.5) -- (-2,2)  
node[anchor=east] {supp $b$};

\draw[densely dotted, thick] (0.5,-0.4) -- (0.8,-2.5) ;
\draw[densely dotted, thick] (-1.45,-0.05) -- (0.8,-2.5) ; 
\draw[densely dotted, thick] (1.45,-0.05) -- (0.8,-2.5) ; 
\draw[densely dotted, thick] (0.3,-1.2) -- (0.8,-2.5)  
node[anchor=west] {$\varpi_h$};

\label{fig1}
\end{tikzpicture}


\BED 
\label{def5.3}
We first fix $R_0 \in \big(\frac14 , \frac12 \big)$ and define for every $h \in (0, \sfh]$, $\sfh = 1/10$,   
the  periodic cell $\varpi_h$ as the union of  the sets
\beas
\varpi_0 &=& \bbD_{R_0}= \{ z \, : \, |z| < R_0  \}  \ \ \ \mbox{and} 
\nonumber \\ 
S_h &=& \{ z \, : \, |{\rm Re}\, z | < 1/2 , \  -h < {\rm Im } \, z < h \} .
\eeas
The periodic domain $\Pi_h$ is defined as the interior of the union
\beas
\bigcup_{m \in \bbZ} \overline{\varpi_h + m }.
\eeas
The Toeplitz operator $T_{\sfa}: A^2(\Pi_h) \to A^2(\Pi_h)$ is determined
by the 1-periodic symbol $\sfa \in L^\infty(\Pi_h)$ (independent of  $h$), 
\begin{equation}
\sfa(z) = \left\{
\begin{array}{ll}
b\big( ( z- [z])/R_0 \big), \ \ \ & z - [z] \in \varpi_0, 
\\
0  & \mbox{otherwise} .
\end{array}
\right.  \label{9.11a}
\end{equation}
For every $h \in (0, \sfh]$ and  $\eta \in [-\pi,\pi]$ we also denote by 
$A_{h,\eta}^2(\varpi_h)$ the space $A^2_\eta(\varpi)$ of Theorem \ref{prop1.2} 
for the periodic domain $\Pi_h$ with periodic cell $\varpi_h$. Also, $P_{h,\eta}$ 
denotes the orthogonal projection from $L^2(\varpi_h)$ onto $A_{h,\eta}^2(\varpi_h)$, 
cf.  Definition \ref{defPe}. We denote the orthogonal projection $P_{\varpi_0}$ from 
$L^2(\varpi_0)$  onto $A^2(\varpi_0)$ briefly by $P_0$. 
\END

The basic idea of the proof of Theorem \ref{9.1} is to use 
Theorem \ref{th8.3} for the domain $\Pi_h$ with a sufficiently small $h$. 
We think of $\varpi_h$ as a domain perturbation of $\varpi_0$, and accordingly
try to approximate the spectra corresponding to $\sigma(T_{a,\eta})$ of \ef{8.40}  
in $\varpi_h$ 
by the spectrum of the restriction of $T_\sfa$ to $\varpi_0$, i.e., the eigenvalues $\lambda_n$ of \ef{9.1} (see also the beginning of Section \ref{secY} for further notation).
The main complication is caused by the fact that although the area of $\varpi_h 
\setminus \varpi_0$ is small, the domains $\varpi_h$ and $\varpi_0$ are of a different character
since the boundary of $\varpi_h$ is non-smooth and it might happen that the 
eigenfunctions in $A_{h,\eta}^2(\varpi_h)$ might be concentrated 
near the corner points of the boundary of $\varpi_h$ and thus they could be 
very different from the eigenfunctions in the smooth domain $\varpi_0$.
These problems can however be circumvented with the help of Lemma \ref{lem9.6},
as will be shown in Section \ref{secY}.

The rest of this section is 
devoted to some necessary facts concerning the related Bergman spaces. The proof of 
Theorem \ref{th9.1} as well as the choice of $h$ will be  completed in the next section.  
If $0 \leq h < k \leq \sfh$, the restriction of functions to a smaller domain
induces  the natural inclusions 
$L^2(\varpi_k) \subset L^2(\varpi_h)$ as well as $A^2(\varpi_k) \subset A^2(\varpi_h)$ and $A_{k,\eta}^2(\varpi_k)
\subset A_{h,\eta}^2(\varpi_h)$, which will be repeatedly used in the sequel.

We will need the following approximation lemma. The difference to standard
approximation results is that the approximating functions need 
to satisfy the quasiperiodic boundary condition \ef{1.15}.

\BEL
\label{lem9.4}
Given $\eta \in [-\pi,\pi]$, the restrictions of the functions 
$f \in A_{\sfh , \eta}^2(\varpi_{\sfh} ) $ to $\varpi_0$ form a dense 
subspace of $A^2(\varpi_0)$. In particular, a given $g \in A^2(\varpi_0)$
can be approximated uniformly on any compact subset of $\varpi_0 $ 
by functions in $A_{\sfh , \eta}^2(\varpi_{\sfh} ) $.
\ENL

Note that the result also holds for all $h \in (0, \sfh]$, since
$ A_{\sfh , \eta}^2(\varpi_{\sfh} ) \subset  A_{h , \eta}^2(\varpi_{h} )$.

\bigskip

Proof. Let $\eta \in [-\pi, \pi]$ and $g \in A^2(\varpi_0)$ be  given. 
The mapping  $\varphi(z) = e^{i 2 \pi  z - 2\pi } $ maps the rectangle 
$Q:=  (-1/2, 1/2) \times (-1/2,1/2) \supset \varpi_\sfh$ conformally onto a subdomain 
$\Omega$ of $\bbD$. We denote the inverse by $\psi : \Omega \to Q$. Moreover,
$\varpi_0 $ is contained in a compact subset of $Q$,
thus, $D:= \varphi ( \varpi_0) $ is contained in a compact subset of 
$\Omega$.

Following Lemma \ref{lem1.5}, one observes that the weighted composition operator 
$J: f \mapsto \varphi' f \circ \varphi $ is a Hilbert-space isomorphism from 
$A^2( D) $ onto $A^2(\varpi_0)$, and $L : A^2(\varpi_0) \to A^2( D) $ with $f \mapsto 
 \psi' f\circ \psi$ is its inverse. Second, $D$ is a Caratheorody domain so that polynomials 
form a dense subspace of $A^2(D)$. 
Finally, it is obvious that $\Omega$ does not contain the line segment
$\{ z \in \bbD \, : \, {\rm Re}\, z < 0 , \ {\rm Im }\, z = 0 \}$, hence,
the logarithm and all real powers of $z$ are defined as analytic functions
on $\Omega$. 

Moreover, the pointwise multipliers $ f \mapsto z^{\eta/(2 \pi)} f$
and  $ f \mapsto \psi' f$ are topological isomorphisms of $A^2(D)$ onto itself:
for some constant
$C>0$ and  all $z \in D$ there holds  $C^{-1} \leq |z|  \leq C$  (since we have 
$D \cap B(0, \rho) = \emptyset$ for a small enough $\rho>0$) and
$C^{-1} \leq |\psi'(z)|  \leq C$ (since the set $D$ is contained in a compact 
subset of the domain of the definition $\Omega$ of the conformal mapping $\psi$).

Consequently, the function  $z^{-\eta/(2 \pi)} g \circ\psi  = z^{-\eta/(2 \pi)} 
L g / \psi'$ belongs to $ A^2(D)$ and we can approximate it by a  polynomial $P$  
with respect to the norm of 
$A^2(D)$. Thus, also $z^{\eta/(2 \pi)} P \psi'$ approximates $ L g $ in $A^2(D)$. Then, due 
to the  isomorphism property of the operator $J$, the function 
\bea
f(z) = J \big( z^{\eta/(2 \pi)} \psi' P  \big)(z) = e^{i \eta z - \eta} P(\varphi(z) )
\in A^2 (\varpi_\sfh)
\eea
approximates $g$ in $A^2(\varpi_0)$, but it is obvious that the entire function
$f$ satisfies the quasiperiodic boundary condition \ef{1.15} and thus also belongs to $A_{\sfh,\eta}^2 (\varpi_\sfh)$.

The uniform approximation in compact subsets is a standard consequence of 
the proven fact. \ \ $\Box$

\bigskip

We prove the following  convergence result for the projection family
$P_{h,\eta}$. Recall that we have fixed $\sfh = 1/10$.

\BEL \label{lem9.6}
Let  $f \in L^2(\varpi_\sfh)$.
Then, for all $\epsilon > 0$ there exists $h_\epsilon \in (0, \sfh]$
such that for all $h \leq h_\epsilon$ and $\eta \in [-\pi,\pi]$ there exist 
$ \tilde g \in A_{h,\eta}^2(\varpi_h)$ with 
\bea
\Vert P_{h,\eta} f -\tilde g \Vert_{\varpi_h} \leq \epsilon
\ \ \ \mbox{and} \ \ \ \Vert P_0f   -\tilde g \Vert_{\varpi_0} \leq \epsilon.
\label{9.25}
\eea
In particular, for every $\eta \in [-\pi,\pi]$ we have $\Vert P_{h,\eta} f - P_0 f\Vert_{\varpi_0} \to 0$ as $h \to 0$.
\ENL

Proof. We may assume $\Vert f \Vert_{\varpi_\sfh} = 1$. Given $\epsilon > 0$,
we find, by  Lemma \ref{lem9.4}, a function $\tilde g \in A_{\sfh,\eta}^2(\varpi_\sfh)$
such that 
\beas
\Vert P_0 f -\tilde g \Vert_{\varpi_0} \leq \epsilon, 
\eeas
which in particular implies the latter inequality in \ef{9.25} and also
that 
\bea
& & \Big( \int\limits_{\varpi_0} | \tilde g - f|^2 dA \Big)^{1/2} 
\leq \Big( \int\limits_{\varpi_0} | \tilde g - P_0f|^2 dA\Big)^{1/2}
+ \Big( \int\limits_{\varpi_0} | P_0f - f|^2 dA\Big)^{1/2}
\rowleq 
\Big( \int\limits_{\varpi_0} | P_0 f - f|^2 dA \Big)^{1/2}+ \epsilon .  \label{9.44}
\eea
Then, we fix $h_\epsilon > 0$ so small  that 
\bea
\int\limits_{\varpi_{h_\epsilon} \setminus \varpi_0} |f - \tilde g|^2 dA \leq \epsilon^2
\ \ \Rightarrow \ \ 
& & \int\limits_{\varpi_h} | \tilde g  - f|^2 dA
\leq  \int\limits_{\varpi_0} | \tilde g - f|^2 dA
+ \epsilon^2  \ \ \forall\ h \leq h_\epsilon . \label{9.42}
\eea
On the other hand, for all $h \leq h_\epsilon$ there holds
\bea
& & \int\limits_{\varpi_h} | P_{h,\eta}f - f|^2 dA
\geq \int\limits_{\varpi_0} | P_{h,\eta}f - f|^2 dA
\geq  \int\limits_{\varpi_0} | P_0f - f|^2 dA  \label{9.46}
\eea
since by the definition of the orthogonal projection, $P_0f$ is the nearest
element of $A^2(\varpi_0)$ to $f$. 
Combining \ef{9.44}--\ef{9.46} yields
\beas
& & \int\limits_{\varpi_h} | \tilde g  - f|^2 dA
\leq \int\limits_{\varpi_h} | P_{h,\eta}f  - f|^2 dA
+ C \epsilon .
\eeas
Thus, the first inequality in \ef{9.25} follows by applying Remark \ref{lem9.3} with 
$H, X, g_1, g_2$ replaced by $L^2(\varpi_h), A_{h,\eta}^2 (\varpi_h), 
P_{h,\eta}f,  \tilde g$, respectively, and changing  $\epsilon$.  \ \ $\Box$

\section{Toeplitz operator with prescribed essential spectrum:
proof.}  \label{secY}

We proceed with the proof of Theorem \ref{th9.1}. Let the symbol $b \in L^\infty(\bbD)$ and $R \in (0,1)$ be given as in the assumptions. We denote by  $\sfb \in 
L^\infty( \varpi_0)$ the function $\sfb(z) = b(z /R_0)$ and also its zero extension to $\varpi_h$. Then, we have 
supp\,$\sfb \subset \bbD_{R R_0}$ so that the support of $\sfb$ is  contained in a compact 
subset of  $\bbD_{R_0} =  \varpi_0 \subset \varpi_h$; see Definition \ref{def5.3}.
Since $z \mapsto z/R_0$ is the conformal mapping from $ \varpi_0$ 
onto $\bbD$, 
Lemma \ref{lem1.5} shows that the eigenvalues of the Toeplitz operator
\bea
T_\sfb : A^2(\varpi_0) \to A^2 (\varpi_0) , \ \ \ \
T_\sfb f = P_0 ( \sfb f )    \label{9.08}
\eea
are the same as the eigenvalues of $T_b: A^2(\bbD) \to A^2(\bbD)$, i.e., the 
numbers in \ef{9.1}. Recall $P_0$ is the Bergman projection of the domain
$ \varpi_0$.

Next, following the notation introduced 
in the previous section, we  denote for every $h \in (0,\sfh]$ and $\eta
\in [-\pi,\pi]$  by $T_{\sfb,h,\eta} : 
A_{h,\eta}^2(\varpi_h) \to A_{h,\eta}^2(\varpi_h)$ the Toeplitz-type operator
\bea 
T_{\sfb,h,\eta} f = P_{h,\eta} (\sfb f) .  \label{9.09}
\eea
We remark that  $T_{\sfb,h,\eta}$ is a compact, self-adjoint operator. The self-adjointness
follows from the property that $\sfb$ is real valued. Moreover, since the support 
of $\sfb$ is compact in $\varpi_h$, a standard normal family argument shows that
the pointwise multiplier $f \mapsto \sfb f$ is a compact operator  from $A^2(\varpi_h)$ 
into $L^2(\varpi_h)$, thus also  from $A_{h,\eta} ^2(\varpi_h)$ into $L^2(\varpi_h)$.
This implies the compactness of $T_{\sfb,h,\eta}$.  Consequently,
we can write the spectrum of $T_{\sfb,h,\eta} : A_{h,\eta}^2(\varpi_h) \to A_{h,\eta}^2(\varpi_h)$
in the same way as in \ef{9.1}, namely, as a union of $\{0\}$ and sequence of eigenvalues 
indexed by taking into account the multiplicities such that 
\bea
|\lambda_{h,\eta}^1| \geq | \lambda_{h,\eta}^2| 
\geq |\lambda_{h,\eta}^3| \geq \ldots \ .
\label{9.11}
\eea

We now proceed via some steps and denote in the sequel by  $\cH$, $\widetilde \cH$ or $
\cH_1,  \cH_2 , \ldots $ infinite subsets of the interval $(0, \sfh]$, all of which have 0 
as an accumulation point.

\smallskip

{\bf Step $1^\circ$.} We fix $\eta \in [-\pi,\pi]$ for the duration of this step. We 
show that if $\lambda \in \bbR $ is an eigenvalue of  $T_b: A^2(\bbD) \to A^2(\bbD)$ 
or equivalently $T_\sfb: A^2(\varpi_0) \to A^2(\varpi_0)$
and $\cH$ is given, then 
for all $h \in  \cH$,
the operator $T_{\sfb,h,\eta}$ has an eigenvalue $\lambda_{h,\eta}$ such that
\bea
\lim_{\stackrel{\scriptstyle h \in \cH,}{  h \to 0}}  \lambda_{h,\eta}
=  \lambda .  \label{9.12}
\eea 
To this end, we only need to consider sufficiently small $h$;  let $0 < \epsilon < 1$ be arbitrary and denote by $f_0 \in A^2(\varpi_0)$, 
$\Vert f_0  \Vert_{\varpi_0} = 1$, an eigenvector  of $T_\sfb$ corresponding to $\lambda$. 
We next use Lemma \ref{lem9.4} to find $\tilde f\in A_{\sfh,\eta}^2 $
such that 
\bea
\Vert \tilde f - f_0 \Vert_{\varpi_0} \leq \epsilon  \label{9.14}
\eea
and then  Lemma \ref{lem9.6} to pick up $h_\epsilon \leq \sfh$ and a function 
$\tilde g \in A_{\sfh,\eta}^2 $ satisfying \ef{9.25} with $f:= \sfb \tilde f \in 
L^2(\varpi_\sfh)$, for all $h \leq h_\epsilon$ and $\eta$. 
By possibly diminishing $h_\epsilon$ we may assume that 
\bea
\Vert \tilde g \Vert_{\varpi_h \setminus \varpi_0} \leq\epsilon
\ \ \mbox{and} \ \ 
\Vert \tilde f \Vert_{\varpi_h \setminus  \varpi_0} \leq\epsilon
\label{9.16}
\eea
for all $h \in \cH$ with $h \leq h_\epsilon$.

 We claim that for these $h$, 
\bea
\Vert T_{\sfb,h,\eta} \tilde f - \lambda \tilde f \Vert_{\varpi_h}
\leq C \epsilon . \label{9.18}
\eea
Note that by \ef{9.14}, \ef{9.16} we have $\big| \, \Vert \tilde f \Vert_{\varpi_h} 
- 1 \big|\leq C \epsilon$, hence, \ef{9.18} and Lemma \ref{lem9.2}  imply that $T_{\sfb,h,\eta}$ has an eigenvalue 
$\lambda_{h,\eta}$ such that $|\lambda_{h,\eta} - \lambda| \leq 
C \epsilon$, which yields the claim \ef{9.12}. 

Now,  \ef{9.16} and \ef{9.25} imply (note $|\lambda| \leq \Vert T_b\Vert $ which is
a fixed constant)
\beas
& & \Vert  T_{\sfb,h,\eta} \tilde f - \lambda \tilde f \Vert_{\varpi_h\setminus \varpi_0}
\leq  \Vert  P_{h, \eta} (\sfb \tilde f)\Vert_{\varpi_h\setminus \varpi_0} 
+ \Vert \lambda  \tilde f \Vert_{\varpi_h\setminus \varpi_0}^2
\rowleq
\Vert  P_{h, \eta} (\sfb \tilde f) -  \tilde g \Vert_{\varpi_h\setminus \varpi_0}^2
+ \Vert  \tilde g \Vert_{\varpi_h\setminus \varpi_0} + C\epsilon  \leq 
C \epsilon.
\eeas
We use this and again \ef{9.25} to obtain $\Vert  P_{h, \eta} (\sfb \tilde f)-
 P_0 (\sfb \tilde f) \Vert_{\varpi_0} \leq C \epsilon$ 
and thus 
\bea
& & \Vert  T_{\sfb,h,\eta} \tilde f - \lambda \tilde f \Vert_{\varpi_h}
\leq  \Vert  P_{h, \eta} (\sfb \tilde f) - \lambda \tilde f \Vert_{\varpi_0}
+ C \epsilon
\rowleq
\Vert  P_0 (\sfb  \tilde f ) - \lambda \tilde f \Vert_{\varpi_0} + C' \epsilon .
\label{9.20}
\eea
Here, \ef{9.14}, the boundedness of the projection operator  and $f_0$ 
being an eigenvector of $T_\sfb$ yield
\beas
\Vert  P_0 (\sfb  \tilde f ) - \lambda \tilde f \Vert_{\varpi_0}
\leq \Vert  P_0 (\sfb  f_0 ) - \lambda f_0 \Vert_{\varpi_0} + C \epsilon
= C \varepsilon,
\eeas
and \ef{9.18} thus follows by taking into account \ef{9.20}.

\smallskip

{\bf Step $2^\circ$.}  
We assume that $\cH$ is a given family and that  $m(h) \in \bbN$, $\eta(h) \in [-\pi,
\pi]$ for all $h \in \cH$. Moreover, it is assumed that there exists 
the limit  
\bea
\lim_{\stackrel{\scriptstyle h \in  \cH,}{  h \to 0}} \lambda^{(h)}= \lambda
\not= 0, 
\ \ \ \mbox{where} \ \lambda^{(h)}:= 
 \lambda_{h,\eta(h)}^{m(h)}, \label{9.30}
\eea
cf. the notation in \ef{9.11}. The  claim is that,  then, $\lambda$ is an eigenvalue of $T_\sfb$.

To see this, we first pick up for each $h \in \cH$  an eigenvector 
$f_h \in  A_{h,\eta(h)}^2$ with $\Vert f_h \Vert_{\varpi_h} = 1$ corresponding to 
the eigenvalue $\lambda^{(h)}$ of the operator $T_{\sfb,h,\eta(h)}$. By a normal family 
argument, there is a decreasing  sequence $\widetilde \cH \subset \cH$  of indices
$h$ converging to zero such that, along this sequence, the functions
$f_h$ converge uniformly on compact subsets of $\varpi_0$ to a function $f_0$, 
which is analytic on $\varpi_0$.  Since $\Vert f_h \Vert_{\varpi_h} = 1$ for all $h$, 
also the limit function satisfies  $\Vert f_0 \Vert_{\varpi_0} \leq  1$. 
Moreover, since the support of $\sfb$ is compact in $\varpi_0$,
we find that $ f_h$ converges uniformly in supp\,$\sfb$ 
to $f_0$. Thus, since $\sfb \in L^\infty(\varpi_0)$, we get
\bea
\sfb f_h \to \sfb f_0 \ \ \ \mbox{in $ L^2(\varpi_0)$ as $h \to 0$.} 
\label{9.30a}
\eea
Now, every projection $P_{h,\eta(h)}$ is a contraction and due to \ef{9.30}, there
exists a constant $C >0$ such that  $\Vert \lambda^{(h)} f_h \Vert_{\varpi_h} 
\geq C$. Thus, we also get for all $h$  a lower bound
\bea
& & \Vert \sfb f_h  \Vert_{\varpi_h}
\geq  \Vert P_{h,\eta(h)} (\sfb f_h ) \Vert_{\varpi_h}
=  \Vert \lambda^{(h)} f_h \Vert_{\varpi_h} 
\geq C > 0  \nonumber \\ 
\Rightarrow  & & \Vert  f_h  \Vert_{\varpi_0}
\geq C' \Vert \sfb f_h  \Vert_{\varpi_0} =  C' \Vert \sfb f_h  \Vert_{\varpi_h} 
\geq C'' > 0.  \label{9.30ab}
\eea
This and \ef{9.30a} imply that $f_0 \not= 0$.

Let $\epsilon > 0$ be arbitrary. We claim that 
\bea
\Vert T_\sfb f_h -  \lambda^{(h)} f_h \Vert_{\varpi_0} \leq \epsilon
\label{9.32}
\eea
for small enough $h \in \widetilde \cH$. If this holds, we pick up $\sigma > 0$
so small that the interval $[ \lambda - \sigma, \lambda + \sigma]$ contains 
at most one eigenvalue of $T_\sfb$. But, in view of \ef{9.30ab}, \ef{9.32} and 
Lemma \ref{lem9.2}, $T_\sfb$ indeed has an eigenvalue
$\tilde \lambda $ with $|\tilde \lambda -  \lambda^{(h)}| \leq C \epsilon$.
Since $\epsilon$ can be arbitrarily small, in particular so that $C \epsilon
< \sigma$, the eigenvalue  $\tilde \lambda$ must be coincide with the limit \ef{9.30}.

To see \ef{9.32}, we again use Lemma \ref{lem9.6} 
to find $h_0> 0$ such that for $h < h_0$ there holds   
\bea
\Vert (P_{h,\eta(h)} - P_0 ) (\sfb f_0 ) \Vert_{\varpi_0} \leq \epsilon ;
\label{9.32a}
\eea
here, $\sfb f_0$ is considered as an element of $L^2(\varpi_h )$
by a zero extension outside $\varpi_0$. We obtain
\beas
& & \Vert T_\sfb f_h   - \lambda^{(h)} f_h \Vert_{\varpi_0} =
\Vert P_0 (\sfb f_h )  - \lambda^{(h)} f_h \Vert_{\varpi_0} 
\rowleq  
\Vert P_0(\sfb f_h) -  P_0 (\sfb f_0  ) \Vert_{\varpi_0} +
\big\Vert ( P_0 - P_{h,\eta(h)} (\sfb f_0 ) ) \big\Vert_{\varpi_0} 
\rowpl
\big\Vert  P_{h,\eta(h)} (\sfb f_0  - \sfb f_h) ) \big\Vert_{\varpi_0}  + 
\Vert  P_{h,\eta(h)} (\sfb f_h  )  - \lambda^{(h)}  f_h \Vert_{\varpi_0} 
. 
\eeas
Here, on the right, the first and third terms are at most $ C\epsilon$ for small enough 
$h$ due to \ef{9.30a}, and for the second one we apply \ef{9.32a}.
The fourth  term is null, since $(\lambda^{(h)}, f_h) $ is an eigenpair of the 
operator $T_{\sfb,h , \eta(h)}= P_{h,\eta(h)} M_\sfb$. This completes step $2^\circ$.

\smallskip

{\bf Step $3^\circ$.} We complete the proof of Theorem \ref{th9.1}. Let $N \in \bbN$,
let $\delta > 0$ be such that $|\lambda_N| > | \lambda_{N+1}| + 2 \delta$
and  let $\varepsilon>0$ be arbitrary. We show that \ef{9.2a} and \ef{9.2b} hold for the operator
$T_\sfa:\Pi_h \to \Pi_h$ with $\sfa$ defined in \ef{9.11a}, if $h > 0$ is small enough.
Indeed, using Step $1^\circ$ (for any fixed $\eta \in [-\pi,\pi]$) and passing 
several times to subsequences, if necessary, we  find a sequence $\cH \subset (0, \sfh]$ 
converging to 0 and, for all  $n=1,\ldots, N$,  indices $m_n(h)$,  such that  
\bea
\lambda_{h,\eta}^{m_n(h)} \to \lambda_n  \label{9.35}
\eea
for all $n=1, \ldots, N$, as $\cH \ni h \to 0$. Now, let $S > 0$ be a  uniform
upper bound for all the operator norms of $T_{\sfb,h,\eta}$, $\eta\in [-\pi,\pi]$,
$h \in (0, \sfh]$, so that the spectra $\sigma(T_{\sfb,h,\eta})$  of all these operators are contained in the 
interval $[-S,S]$. Let also $\epsilon' > 0$ with $\epsilon' < \varepsilon$ be such that all  closed balls  $\overline{B(\lambda_n, \epsilon')}$ 
are disjoint for 
$n =1, \ldots, N$ (except for those where $\lambda_n = \lambda_\ell$ for some $n \not=\ell$,
when $B(\lambda_n, \epsilon') = B(\lambda_\ell, \epsilon'))$.  Note that by the choice of $\delta $, the set
\bea
 \bfS_N:=  \big[-S , -|\lambda_{N+1}| - \delta\big] \cup 
\big[ |\lambda_{N+1}| +  \delta , S\big] \subset \bbR
\label{9.36}
\eea
contains only the eigenvalues $\lambda_n$ with $n=1, \ldots,N$ and no others.
Thus, the set 
$$
\Sigma_0 := \sigma(T_b) \cap \Big( \bfS_N \setminus \bigcup_{n=1}^N  
B(\lambda_n, \epsilon') \Big) 
$$ 
is empty. This leads to the conclusion
that there must exist $h_0 > 0$ such that also the set
\beas
\Sigma_{h,\eta}:= \sigma(T_{\sfb,h,\eta})  \cap 
\Big( \bfS_N \setminus \bigcup_{n=1}^N B(\lambda_n, \epsilon') \Big) 
\eeas
is empty for all $h < h_0$ and all $\eta \in [-\pi,\pi]$. Otherwise we would find a 
sequence of numbers $\lambda_{h_k, \eta_k} \in \Sigma_{h_k,\eta_k} \subset 
\sigma(T_{\sfb,h_k,\eta_k})$ with $h_k \to 0$ as $k \to \infty$, also contained in 
the compact set $\bfS_N \setminus \cup_{n=1}^N  B(\lambda_n, \epsilon')$. Passing to a subsequence if necessary, we could assume that the sequence $\big( \lambda_{h_k, \eta_k} \big)_{k=1}^\infty$ would converge (to a non-zero point, since $\bfS$ 
is bounded away from zero).
According to  Step $2^\circ$, the limit would be an eigenvalue of $T_b$ and thus a point
in the empty set $\Sigma_0$, which is  impossible.

Now, we take a fixed $h > 0$ so small that $h < h_0$ and $\lambda_{h,\eta}^{m_n(h)}
\in B(\lambda_n, \epsilon')$ for all $n = 1, \ldots, N$, see \ef{9.35}.
Then, the domain $\Pi_h$ is chosen for the domain $\Pi$ in Theorem \ref{th9.1},
and the operator $T_\sfa$ is as defined in \ef{9.11a}. 
Due to formula \ef{8.40}, the relation \ef{9.2b} holds. 
Relation \ef{9.2a} follows from the choice of $\delta$, \ef{9.1}, \ef{9.36}  and
the fact that the set $\Sigma_{h,\eta}$ is empty. This  completes the proof of Theorem \ref{th9.1}.
\ \ $\Box$

\end{document}